\documentclass[12pt,one,english]{amsart}
\usepackage{mathabx}
\usepackage{braket}
\usepackage{bbm}
\usepackage[T1]{fontenc}
\usepackage{amsfonts}
\usepackage{mathrsfs}
\usepackage{stmaryrd}
\usepackage{upgreek}
\usepackage{textcomp}
\usepackage[utf8]{inputenc}

\usepackage{extarrow}
\usepackage{enumerate}
\usepackage{enumitem}
\usepackage{graphicx}
\usepackage[colorlinks=true, linkcolor=blue]{hyperref}
\usepackage{latexsym}
\usepackage{mathtools}
\usepackage{mathtext}
\usepackage{amsmath}
\usepackage{amssymb}
\usepackage{amsthm}
\usepackage{tikz}
\usepackage{tikz-cd}
\usepackage[all,cmtip]{xy}
\usepackage[mathscr]{eucal}
\usepackage[toc,page]{appendix}
\usepackage{a4wide}
\usepackage[
backend=biber,
style=alphabetic-verb,
]{biblatex}

\usepackage[a4paper, left=2cm, right=2cm, top=2.5cm, bottom=2.5cm]{geometry}

\usetikzlibrary{arrows}
\usetikzlibrary{matrix}
\usetikzlibrary{shapes}
\usetikzlibrary{snakes}
\usetikzlibrary{matrix}

\DeclareMathOperator{\End}{\textup{End}}

\DeclareMathOperator{\FEt}{\textup{FÉt}}

\DeclareMathOperator{\Sh}{\textup{Sh}}

\DeclareMathOperator{\Spec}{\textup{Spec}}

\DeclareMathOperator{\et}{\textup{et}}

\theoremstyle{plain}
\newtheorem{thm}{Theorem}[section]
\newtheorem*{thm*}{Theorem}
\newtheorem{lem}[thm]{Lemma}
\newtheorem{cor}[thm]{Corollary}
\newtheorem{prop}[thm]{Proposition}

\theoremstyle{remark}
\newtheorem{rmk}[thm]{Remark} %\setcounter{rmk}{0}

\newtheorem*{rmk*}{Remark}

\theoremstyle{definition}
\newtheorem{defn}{Definition}[section] 

\newtheorem*{const*}{Construction}

\newtheorem{setting}[defn]{Setting}

\theoremstyle{plain}
\newtheorem{thmI}{Theorem}
\def \fppf {{\textup{lis-ét}}}
\def \gp {{\rm gp}}
\def \sh {{\mathit{Sh}}}
\def \flc {{\rm flc}}

\AtEveryBibitem{\clearlist{language}}

\global\long\def\F{\mathbb{F}}

\global\long\def\N{\mathbb{N}}

\global\long\def\Z{\mathbb{Z}}

\global\long\def\id{\textup{id}}

\newcommand{\sA}{{\mathcal A}}
\newcommand{\sB}{{\mathcal B}}
\newcommand{\sC}{{\mathcal C}}

\newcommand{\sF}{{\mathcal F}}

\newcommand{\sO}{{\mathcal O}}

\newcommand{\sT}{{\mathcal T}}
\newcommand{\sU}{{\mathcal U}}

\newcommand{\sX}{{\mathcal X}}
\newcommand{\sY}{{\mathcal Y}}
\newcommand{\sZ}{{\mathcal Z}}

\newcommand{\an}{{\rm an}}

\renewcommand{\et}{\textup{\'et}}

\newcommand{\Hom}{{\rm Hom}}

\newcommand{\tp}{{\rm top}}

\begin{document}
\title{Drinfeld's Lemma for Algebraic Stacks}
\author{Lei Zhang}

 \address{Lei Zhang\\
     Sun Yat-Sen University\\ School of Mathematics
     (Zhuhai)\\Zhuhai, Guangdong, P.~R.~China}
\email{cumt559@gmail.com}

%\thanks{This work was supported by the European Research Council (ERC) Advanced Grant 0419744101 and the Einstein Foundation}
\date{\today}

\makeatletter 
\providecommand\@dotsep{5} 
\makeatother 
%\listoftodos\relax

\setcounter{section}{0}
\maketitle

\begin{abstract} Drinfeld's lemma is a powerful tool for splitting $\ell$-adic
    local systems defined over a product of connected schemes over a finite
    field. We extend the lemma to products of connected algebraic stacks when
    all but possibly one factor are quasi-compact and quasi-separated. The
    associated category of finite \'{e}tale covers with partial Frobenius
    actions is a Galois category without any quasi-compactness assumption.
\end{abstract}

\section{Introduction}

The main motivation of Drinfeld's lemma is to split $\ell$-adic local systems
defined over a product of schemes. More precisely, let $X_1,X_2$ be two connected
$\F_q$-schemes, then one would like to get, out of an $\ell$-adic local system
%in $\Loc_X(\bar{\Q}_{\ell})$
on
$X\coloneq X_1\times_{\F_q}X_2$, an $\ell$-adic local system coming from
local systems on the individual factors
$X_1,X_2$. The problem is
easy if one considers
complex analytic local systems on a product of complex varieties. Indeed, one
can split the local system via the Künneth formula for
topological fundamental groups:
\[
    \pi_1^{\tp}(X_1^\an\times X_2^\an)\xlongrightarrow{\cong}\pi_1^{\tp}(X_1^\an)\times\pi_1^{\tp}(X_2^\an)
\]
One can also do this for $\ell$-adic local systems. Indeed, for
any connected qcqs schemes $X_1,X_2$ defined over an algebraically closed field $k$
\emph{of characteristic 0}, one has an isomorphism: 
\[
    \pi_1^{\et}(X_1\times_k
    X_2)\xlongrightarrow{\cong}\pi_1^{\et}(X_1)\times\pi_1^{\et}(X_2)
\]
Even when $k$ is of characteristic $p>0$ (but still algebraically closed), the above Künneth
formula still holds true provided that either $X_1$ or $X_2$ is proper over $k$.
However, the Künneth formula fails when $k$ is not
algebraically closed: take $X_1=X_2=\Spec(\F_p)$, then the Künneth formula
would mean that the diagonal map $\hat{\Z}\to\hat{\Z}\times\hat{\Z}$ is an
isomorphism. 

The issue (for finite fields) can be resolved if partial Frobenii actions are brought into play. More
precisely, let $\phi_1$ (resp. $\phi_2$) denote the \textit{partial Frobenius}
map
$X\to X$, which is the $q$-absolute Frobenius on $X_1$
(resp. $X_2$) and the identity on the other. Consider the category $\FEt(X/\Phi)$ of triples
$(Y,\varphi_1,\varphi_2)$,
where $Y$ is a finite étale cover of $X=X_1\times_{\F_q}X_2$ and $
\varphi_i$ is an isomorphism $Y\xrightarrow{\cong} \phi_i^*Y$ satisfying
that
$\phi^*_1(\varphi_2)\circ\varphi_1=\phi^*_2(\varphi_1)\circ\varphi_2$
is the identity (by identifying $(\phi_1\circ\phi_2)^*Y$ with $Y$ via the
absolute Frobenius of $Y$).

\begin{thm}[Drinfeld's lemma for schemes]\label{Drinfeld's lemma for schemes} Suppose $X_1,X_2$ are connected
    quasi-compact and quasi-separated (qcqs) $\F_q$-schemes.
    Then \begin{itemize}
        \item $\FEt(X/\Phi)$ is a Galois category whose Galois group is denoted
            by $\pi_1^\et(X/\Phi)$;
        \item the  natural map $\pi_1^\et(X/\Phi)\to
            \pi_1^\et(X_1)\times\pi_1^\et(X_2)$ is an isomorphism.
    \end{itemize}
\end{thm}

Similarly, one has Drinfeld's lemma for $n$ factors $X_1,\dots, X_n$. Please refer to \cite[Theorem 2.1]{Dri80}, \cite[IV.2, Theorem 4]{Laf97}, \cite[Theorem 8.1.4]{Lau07},
\cite[Theorem 4.2.12]{Ked17}, \cite[Lemma 0.18]{Laf18}, \cite[Theorem 16.2.4]{SW20},  and \cite[Theorem
1.4]{Müller22} for details. Using Drinfeld's Lemma one can split local systems
on a product of schemes
equipped with commuting partial Frobenius actions.

The notion of fundamental group of algebraic stacks was introduced and
studied by B.~Noohi in \cite{Noohi00}. In this note we generalize
Theorem \ref{Drinfeld's lemma for schemes} to algebraic
stacks and add details to its proof. We also allow one factor to be non-qcqs;
the other factors are required to be qcqs. This qualification is necessary
for the spreading-out argument used below, while the existence of the Galois
category itself requires no qcqs hypothesis.
\begin{thmI}[cf.~\S\ref{sec: general remarks} and \S\ref{Drinfeld's lemma}]\label{Drinfeld's Lemma}
Let $n\geq 1$, let $\sX_1,\sX_2,\dots,\sX_n$ be connected algebraic stacks over $\F_q$, and set
$\sX\coloneqq \sX_1\times_{\F_q}\dots\times_{\F_q}\sX_n$. Then
\begin{enumerate}[label=\rm{(\arabic*)}]
    \item 
The category
\(\FEt(\sX/\Phi)\) consisting  of
finite étale covers of $\sX$ equipped
with commuting partial Frobenii actions is a Galois category, whose Galois group is denoted by
$\pi_1^\et(\sX/\Phi)$;
\item if all but possibly one of the stacks $\sX_i$ are qcqs, the natural map $\pi_1^\et(\sX/\Phi)\longrightarrow
    \pi_1^\et(\sX_1)\times\cdots\times\pi_1^\et(\sX_n)$ is an isomorphism.
\end{enumerate}
\end{thmI}

As in the scheme-theoretic proof of Drinfeld's lemma, the permanence of
$\pi_1^\et$ -- corresponding to the special case of Theorem \ref{Drinfeld's
Lemma} where $n=2$ and $\sX_2 = \Spec(\bar{k})$ for a field $k$ -- plays a central role. Here, we employ its stack-theoretic analogue (Lemma \ref{Permanence}), which, when combined with standard scheme-theoretic arguments and dévissage techniques, yields the desired result.
A key preliminary step is to properly define the category $\FEt(\sX/\Phi)$ for an algebraic stack $\sX$. We formulate this in greater generality:
Let $M$ be a commutative monoid acting on an algebraic stack $\sX$, and let $G\coloneqq M^\gp$ denote its Grothendieck group. Define $\FEt(\sX/M)$ as the category of finite étale covers $\sY \in \FEt(\sX)$ equipped with an $M$-action compatible with the action on $\sX$. When the quotient topological space $|\sX|/G$ is connected, $\FEt(\sX/M)$ forms a Galois category (Theorem \ref{Galois category}).
This framework necessitates working with monoid actions on stacks. While M. Romagny's \cite{Rom05} treats group actions systematically, the same definitions apply verbatim to monoids. Notably, in our setting (Lemma \ref{universal homeomorphism vs cartesian}), the action simplifies considerably: higher compatibility conditions become automatic, bypassing much of the complexity in Romagny's general theory.

\bigskip

\noindent\textbf{Notation and conventions}
\begin{enumerate}[label=(\arabic*)]
    %\item All \emph{algebraic stack} over a base scheme $S$ in this paper are fibered over the fppf-site $(\Sch/S)_{\mathrm{fppf}}$ (cf.~\cite[\href{https://stacks.math.columbia.edu/tag/026O}{026O}]{stacks-project}).
    \item Throughout this paper the word \emph{representable} means
        representable by algebraic spaces.
    \item By a \emph{finite} map $f\colon \sY\to\sX$ of algebraic stacks, we mean a
        map which is representable by a finite map of schemes
        (cf.~\cite[\href{https://stacks.math.columbia.edu/tag/0CHU}{0CHU}]{stacks-project}). 
\end{enumerate}

\section{The $G$-Connectedness}

One of the difficulties in understanding Drinfeld's lemma is that the partial $q$-Frobenius maps are, in general,  not invertible.
Thus ``the quotient space $X/\Phi$'' is only a suggestive symbol but not an
honest space -- even when $X$ is a scheme!
However, we can first forget
this if we only look at the actions on the ambient topological space $|\sX|$.

\begin{defn} Let $G$ be a group, and let $X$ be a topological space equipped
with a $G$-action $\rho_G\colon G\times X\to X$. The pair $(X,\rho_G)$ is called a
\emph{$G$-space}. The space $X$ is called \emph{$G$-connected} if the quotient
space $X/G$ is a connected topological space. A subspace $S\subseteq X$ is
called a \emph{$G$-component} of $X$ if it is the inverse image of a nonempty
clopen subspace of $X/G$. An ordinary
topological space $X$ can be viewed as a $G$-space via the trivial $G$-action.
In this case, the pair $(X,\rho_G)$ will simply be denoted by $X$.\end{defn}

Let $f\colon (X,\rho_G)\to (S,\rho_G')$ be a map, \emph{i.e.} a map of topological spaces
$f\colon X\to S$ which is compatible with the $G$-actions. For each $G$-stable
subset $S_1\subseteq S$ the inverse image $f^{-1}(S_1)\subseteq X$ is a $G$-stable subspace. 
\begin{lem}\label{geometrically phi-connected}
  Let $f\colon (X,\rho_G)\to S$ be map whose fibers are $G$-connected. Suppose
  that $X\to S$ is submersive. Then $(X,\rho_G)$ is $G$-connected if and only if $S$ is connected.
\end{lem}
\begin{proof} The map $f$ induces a continuous map $\bar{f}\colon X/G\to S$
whose fibers are connected (hence nonempty). If $(X,\rho_G)$ is $G$-connected, then $X/G$ is connected,
so $S$ is connected. Conversely, if $(X,\rho_G)$ is not $G$-connected,
then $X/G=U_1\sqcup U_2$, where $U_1,U_2$ are two non-empty
open subspaces of $X/G$. Since the fibers of $\bar{f}$ are
connected, any fiber is contained either in $U_1$ or in $U_2$. Thus exists subspaces $V_1,V_2$ of $S$ such that $S=V_1\sqcup V_2$ and $\bar{f}^{-1}(V_i)=U_i$ ($i=1,2$). That $f$ is submersive implies that $\bar{f}$ is submersive. Hence $V_1,V_2$ are non-empty opens of $S$. Then  $S$ is not connected.
\end{proof}

\section{Frobenius and universal homeomorphisms}

Recall that for any fibered category $\sX$ defined over $\F_p$, 
the absolute Frobenius map $F_{\sX}\colon \sX\to\sX$ sends a
section $x\in\sX(T)$ on an
$\F_p$-scheme $T$ to the composition
$T\xrightarrow{F_T}T\xrightarrow{x}\sX\in\sX(T)$, where
$F_T$ is the absolute Frobenius of the scheme $T$.

Universal homeomorphisms of schemes are radicial
morphisms
(cf.~\cite[\href{https://stacks.math.columbia.edu/tag/01S2}{01S2}]{stacks-project})
which play a role as the ``opposite'' of étale maps. The on the nose extension
of this notion to stacks could be misleading in certain situations.
 For example, consider a  finite étale group scheme
$G$ over a field $k$. Any base change of the map $f\colon\sB_kG\to\Spec(k)$
is a homeomorphism on the underlying topological spaces, however, the map is
 étale. This
suggests a new definition of the notion of ``universal homeomorphism'' for
algebraic stacks -- we should consider not only $f$ but also the diagonal of
$f$ which encodes information about the (relative) inertia.

\begin{defn} A map of algebraic stacks $f\colon \sX\to\sY$ is called a
    \textit{universal homeomorphism on the nose} if for any map of algebraic stacks
    $\sT\to\sY$, the base change $f_{\sT}\colon \sX\times_{\sY}\sT\to\sT$ is a
    homeomorphism on the underlying topological
    spaces. The map $f$ is called a
    \textit{universal homeomorphism} if both $f$ and the
    diagonal $\Delta_f$ are universal homeomorphisms on the nose.
\end{defn}
\begin{rmk} The term ``universal homeomorphism on the nose'' corresponds to
    what is simply called a ``universal homeomorphism'' in
    \cite[\href{https://stacks.math.columbia.edu/tag/0CI7}{0CI7}]{stacks-project}.
    However, as will become apparent in the proof of Lemma \ref{decent
    geometric properties of uni homeomorphisms} (2) and Lemma \ref{universal
    homeomorphism vs cartesian}
    (a)$\Leftrightarrow$(b),
    this standard notion proves inadequate for our purposes.
\end{rmk}

\begin{lem} If $\sX$ is an algebraic stack over $\F_p$, then $F_{\sX}$ is a
    universal homeomorphism.
\end{lem}
\begin{proof}
    Note that for any perfect field $k$, the absolute Frobenius $F_{\sX}$ induces an
    equivalence:
    $F_{\sX}(k)\colon\sX(k)\xrightarrow{\cong}\sX(k)$, and its (strict) inverse is given
    by sending $\Spec(k)\xrightarrow{x}\sX$ to
    $\Spec(k)\xrightarrow{F_k^{-1}}\Spec(k)\xrightarrow{x}\sX$. As an immediate
    consequence $F_{\sX}$ induces a bijection $|F_{\sX}|\colon|\sX|\xrightarrow{\cong}|\sX|$.
Suppose $f\colon \sY\to\sX$ is a map of algebraic stacks over $\F_p$. Then we
have the following cartesian diagram
\[
    \begin{tikzpicture}[xscale=2.0,yscale=1.2,bmr/.pic={\draw (0,0)--++(-90:2mm)--++(180:2mm);}]
         \path
                (0,0)     node (F) {$\sY'$}
             +(-43:.4) pic[scale=1,red]{bmr}
             +(0:1.5)  node (star) {$\sY$}
             ++(-90:1.5) node (X) {$\sX$}
             +(0:1.5)  node (Y) {$\sX$};
         \draw[->] (F)--(star) node[midway,above,scale=.6]{$F'$};
         \draw[->] (F)--(X) node[midway,left,scale=.6]{$f'$};
         \draw[->] (X)--(Y) node[midway,above,scale=.6]{$F_{\sX}$};
         \draw[->] (star)--(Y) node[midway,right,scale=.6]{$f$};
     \end{tikzpicture}
\] 
The universality of a 2-fibered square implies that $\sY'(k)\simeq
\sX(k)\times_{F_{\sX}(k)}\sY(k)$ for any perfect field $k$. Thus $F'(k)\colon\sY'(k)\to\sY(k)$ is also an
equivalence.
This shows that $F'$ induces a bijection $|F'|\colon|\sY'|\to |\sY|$.

Let $F_{\sY/\sX}\colon \sY\to\sY'$ be the relative Frobenius map. 
We have seen that $|F'|$, $|F_{\sY/\sX}|$, $|F_{\sY}|$ are bijective
and continuous. From the very definition of the topology on $|\sY|$ it is clear
that $|F_{\sY}|$
is a homeomorphism. The continuity of $|F_{\sY/\sX}|$ then implies that $|F'|$ is
submersive, hence a homeomorphism.

Let's consider the diagonal $\Delta\colon \sX\to
\sX\times_{F_{\sX}}\sX$. Again, for any perfect field $k$, $\Delta(k)$ equals
the diagonal map $\sX(k)\to\sX(k)\times_{F_{\sX}(k)}\sX(k)$ which is an
equivalence.  Now suppose we have a
pullback diagram 
\[
    \begin{tikzpicture}[xscale=2.0,yscale=1.2,bmr/.pic={\draw (0,0)--++(-90:2mm)--++(180:2mm);}]
         \path
                (0,0)     node (F) {$\sY'$}
             +(-43:.4) pic[scale=1,red]{bmr}
             +(0:1.5)  node (star) {$\sY$}
             ++(-90:1.5) node (X) {$\sX$}
             +(0:1.5)  node (Y) {$\sX\times_{F_{\sX}}\sX$};
         \draw[->] (F)--(star) node[midway,above,scale=.6]{$\Theta$};
         \draw[->] (F)--(X) node[midway,left,scale=.6]{$$};
         \draw[->] (X)--(Y) node[midway,above,scale=.6]{$\Delta$};
         \draw[->] (star)--(Y) node[midway,right,scale=.6]{$f_1\times f_2$};
     \end{tikzpicture}
\] 
Then again, $\Theta(k)\colon \sY'(k)\to\sY(k)$ is an equivalence. Thus
$|\Theta|$ is bijective. To show that $|\Theta|$ is a homeomorphism, one
notices the following 2-commutative diagram.
\[
    \begin{tikzpicture}[xscale=2.0,yscale=1.2,bmr/.pic={\draw (0,0)--++(-90:2mm)--++(180:2mm);}]
         \path
                (0,0)     node (F) {$\sY$}
             %+(-43:.4) pic[scale=1,red]{bmr}
             +(0:1.5)  node (star) {$\sY$}
             ++(-90:1.5) node (X) {$\sX$}
             +(0:1.5)  node (Y) {$\sX\times_{F_{\sX}}\sX$};
             \draw[->] (F)--(star) node[midway,above,scale=.6]{$F_{\sY}$};
             \draw[->] (F)--(X) node[midway,left,scale=.6]{$f_1\circ F_{\sY}$};
         \draw[->] (X)--(Y) node[midway,above,scale=.6]{$\Delta$};
         \draw[->] (star)--(Y) node[midway,right,scale=.6]{$f_1\times f_2$};
     \end{tikzpicture}
 \] The homeomorphism $|F_{\sY}|$ factors through the bijective continuous
 map $|\Theta|$, then $|\Theta|$ must be a homeomorphism as well.
\end{proof}

Here are some geometric properties of universal homeomorphisms for algebraic
stacks. Note that Lemma \ref{decent geometric properties of uni
homeomorphisms} (2) is not true if $g$ was only a universal
homeomorphism on the nose.

\begin{lem} \label{decent geometric properties of uni homeomorphisms}
    Let $\sX\xrightarrow{f}\sY\xrightarrow{g}\sZ$ be a sequence of maps of
    algebraic stacks. Set $h\coloneqq g\circ f$. 
    \begin{enumerate}[label={\rm(\arabic*)}]
         \item If $f,g$ are universal homeomorphisms, then so is $h$.
         \item If $g$ is a universal homeomorphism and  $h$
             is a universal homeomorphism on the nose, then $f$ is a
            universal homeomorphism on the nose.
        \item If $\sX,\sY$ are algebraic spaces, then $f$ is a universal
            homeomorphism iff it is a universal homeomorphism on the nose.
        % \item If $g,h$ are universal homeomorphisms, then $f$ is a universal homeomorphism.
    \end{enumerate}
\end{lem}
\begin{proof} (1), (3) are  obvious. Let's look at (2). The map $f$ factorizes
    as
$\sX\xrightarrow{\Gamma_f}\sX\times_{\sZ}\sY\xrightarrow{h\times\id_{\sY}}\sY$,
    where $\Gamma_f$ is the graph of $f$. If $h$ is a universal homeomorphism
    on the nose and  $g$ is a universal homeomorphism, then  $\Gamma_f$
    and $h\times\id_{\sY}$ are universal homeomorphisms on the nose, and so is
    the composition $f$.
\end{proof}

Note that the absolute Frobenius is only $\F_p$-linear. In order to make it
$\F_q$-linear we take the $q$-absolute Frobenius map $\phi_{\sX}\coloneqq F_{\sX}^n$,
where $n\coloneqq \log_pq$. Just as in the scheme case, the pullback of the étale
maps along the $q$-absolute Frobenius
map of algebraic stacks induces an equivalence of étale sites (actually the
identity!).

\begin{prop} Let $\sX$ be a Deligne-Mumford stack {\rm(}resp. an algebraic
    stack{\rm)}
    over $\F_q$, and let $f\colon \sY\to \sX$ be an étale map {\rm(}resp. a representable
    étale map{\rm)}. Then we have the following cartesian diagram.
    \[
            \begin{tikzpicture}[xscale=2.0,yscale=1.2,bmr/.pic={\draw (0,0)--++(-90:2mm)--++(180:2mm);}]
                 \path
                        (0,0)     node (F) {$\sY$}
                     +(-43:.4) pic[scale=1,red]{bmr}
                     +(0:1.5)  node (star) {$\sY$}
                     ++(-90:1.5) node (X) {$\sX$}
                     +(0:1.5)  node (Y) {$\sX$};
                 \draw[->] (F)--(star)
                     node[midway,above,scale=.6]{$\phi_{\sY}$};
                 \draw[->] (F)--(X) node[midway,left,scale=.6]{$f$};
                 \draw[->] (X)--(Y) node[midway,above,scale=.6]{$\phi_{\sX}$};
                 \draw[->] (star)--(Y) node[midway,right,scale=.6]{$f$};
            \end{tikzpicture} 
    \]
\end{prop}
\begin{proof} One just has to show that the natural map $\lambda\colon\sY\to
    \sX\times_{\phi_{\sX}}\sY$ is an equivalence. For this, one can take an
    atlas of $X\twoheadrightarrow \sX$, where $X$ is a scheme and show that the
    pullback of $\lambda$ to $X$ is an equivalence. Thus we are reduced to the
    case when $\sX=X$ is a scheme. If $f$ is representable by schemes, then we
    are done. If $f$ is only representable by algebraic spaces, then $\sY=Y$ is
    an algebraic space. Take an étale atlas $Z\twoheadrightarrow Y$. The following diagram
    \[ 
        \begin{tikzpicture}[xscale=2.0,yscale=1.2,bmr/.pic={\draw (0,0)--++(-90:2mm)--++(180:2mm);}]
             \path
                    (0,0)     node (F) {$Z$}
                 +(-43:.4) pic[scale=1,red]{bmr}
                 +(0:1.5)  node (star) {$X\times_{\phi_X}Z$}
                 ++(-90:1.5) node (X) {$Y$}
                 +(0:1.5)  node (Y) {$X\times_{\phi_X}Y$};
                 \draw[->] (F)--(star)node[midway,above,scale=.6]{$\lambda_{Z}$};
             \draw[->>] (F)--(X)node[midway,above,scale=.6]{};
             \draw[->] (X)--(Y) node[midway,above,scale=.6]{$\lambda$};
             \draw[->>] (star)--(Y)node[midway,above,scale=.6]{};
        \end{tikzpicture} 
    \] is cartesian, because its extension via the projections
    $X\times_{\phi_X}Z\to Z$ and $X\times_{\phi_X}Y\to Y$ is cartesian (as
    $Z\twoheadrightarrow Y$ is representable by schemes). Since $\lambda_Z$ is
    an isomorphism by the scheme case, $\lambda$ is an isomorphism as well. If
    $\sY$ is only a Deligne-Mumford stack, then we take an étale atlas
    $Z\twoheadrightarrow \sY$ and repeat the above argument to reduce the
    problem to the case when $f$ is representable by algebraic spaces which has
    just been proven.
\end{proof}

\section{The $G$-connectedness under partial Frobenius maps}
Let's come back to our original setting. Let $\sX_1,\sX_2,\dots,\sX_n$ be connected algebraic stacks over $\F_q$, and set $\sX\coloneqq \sX_1\times_{\F_q}\cdots\times_{\F_q}\sX_n$.
Then there are \textit{partial Frobenius
endomorphisms}  $\phi_i\colon \sX\to \sX$, which is the $q$-absolute Frobenius
on $\sX_i$
 and the identity on the others. These $\phi_i$s are commuting endomorphisms of
 $\sX$ (resp. automorphisms of $|\sX|$) whose composition
 $\phi_1\circ\cdots\circ\phi_n$ is the $q$-Frobenius of $\sX$ (resp. the
 identity of $|\sX|$). 

 In some cases, e.g. when studying the actions of the $\phi_i$s on $|\sX|$ or
 $\FEt(\sX)$, one
 can drop any of the partial Frobenius maps (say $\phi_{n}$), as $\phi_n$ is completely determined by all the
 others. The actions of $\{\phi_i\}_{1\leq i\leq n-1}$ provide  $\sX$ (resp. $|\sX|$)  with a strict
 action (resp. an action) $\rho_{\sX}$ by the monoid
 $M_0\coloneqq\phi_1^\N\times\cdots\phi_{n-1}^\N=\N^{n-1}$ (resp. the
 group $G_0\coloneqq\phi_1^\Z\times\cdots\phi_{n-1}^\Z=\Z^{n-1}$).

\begin{cor}\label{connectedness}
The space $|\sX_1\times_{\F_q}\cdots\times_{\F_q}\sX_n|$
equipped with the $\Z^{n-1}$-action by $n-1$ of the partial Frobenius maps is
$\Z^{n-1}$-connected. 
\end{cor}
\begin{proof} Consider the projection map
    $|\textup{pr}_n|\colon|\sX_1\times_{\F_q}\cdots\times_{\F_q}\sX_n|\to
    |\sX_n|$, where the former is equipped with the $\Z^{n-1}$-action via the $\phi_{i}$s,
    and the latter is equipped with the trivial $\Z^{n-1}$-action. Let $x\in\sX_n(k)$ be a
    geometric point representing a point in $|\sX_n|$. Then the topological space of the
    stack fiber $\sX_1\times_{\F_q}\cdots\times_{\F_q}\sX_{n-1}\times_{\F_q}k$ maps surjectively to the topological fiber of $x$ along the map $|\textup{pr}_n|$. This surjection moreover preserves the $\Z^{n-1}$-actions. 
    In light of Lemma \ref{geometrically phi-connected}, it is enough to show that the space $|\sX_1\times_{\F_q}\cdots\times_{\F_q}\sX_{n-1}\times_{\F_q}k|$ is $\Z^{n-1}$-connected. 

    Suppose $\sU\subseteq \sX_1\times_{\F_q}\cdots\times_{\F_q}\sX_{n-1}\times_{\F_q}k$
    is a clopen substack which is stable under the $\Z^{n-1}$-action.
    Then $\sU$ is also stable under the $\Z^{n-1}$-action defined by
    $\phi_1,\dots,\phi_{n-2},\phi_k$. By \cite[Theorem 4.1]{DTZ20}, $\sU$
    is the preimage of a clopen substack $\sU_0\subseteq \sX_1\times_{\F_q}\cdots\times_{\F_q}\sX_{n-1}$,
    and $\sU_0$ is obviously stable under the $\Z^{n-2}$-action
    $\phi_1,\dots,\phi_{n-2}$.
    Then $\sU_0$ is either empty or the whole stack by the connectedness of $\sX_{n-1}$ in
    the case $n=2$, or by the induction hypothesis in the case $n>2$. Thus $\sU$ is also either empty
    or the whole stack. This finishes the proof.
\end{proof}

\section{Galois categories}
\label{sec: general remarks}

Let $\sX$ be an algebraic stack over $\F_q$. When $\sX$ is connected,
it is known
(cf.~\cite[Theorem 4.2]{Noohi00} or \cite[\S 4]{Noohi04}) that $\FEt(\sX)$ is a
Galois category.
Let's now adapt this to our setting. Suppose $M$ is a commutative monoid
acting (cf.~\cite[Def.~1.3 (i)]{Rom05}) on $\sX$ via universal
homeomorphisms. Let $\rho_{\sX}\colon M\to\End_{\F_q}(\sX)$ be the action map. If $G$ denotes the Grothendieck group associated with
$M$, \emph{i.e.} $G\coloneqq M^\gp$,
then $|\sX|$ is a $G$-space.
\begin{defn} Let $\FEt(\sX/M)$ denote the category whose
objects are 2-commutative
diagrams
\[
            \begin{tikzpicture}[xscale=1.5,yscale=0.9]
                 \path
                        (0,0)     node (F) {$\sY$}
                     
                     +(0:1.5)  node (star) {$\sY$}
                     ++(-90:1.5) node (X) {$\sX$}
                     +(0:1.5)  node (Y) {$\sX$};
                     \node (start) at (1.1,-0.5) {};
                     \node (end) at (.5,-1.1) {};
                 \draw[->] (F)--(star) node[midway,above,scale=.6]{$\rho_{\sY}(m)$};
                 \draw[->] (F)--(X) node[midway,left,scale=.6]{$f$};
                 \draw[->] (X)--(Y) node[midway,above,scale=.6]{$\rho_{\sX}(m)$};
                 \draw[->] (star)--(Y) node[midway,right,scale=.6]{$f$};
                 \draw[double distance=3.5pt,line width=.8pt,-{Classical TikZ
                 Rightarrow[length=3mm, line width=.8pt]}] (start)--(end)
                 node[midway,above,scale=.6,yshift=7pt]{$\sigma_m$};
            \end{tikzpicture}
        \]
        where $\sY\xrightarrow{f}\sX\in\FEt(\sX)$ and
$\rho_{\sY}$ is an 
    action of $M$ on $\sY$ via $\F_q$-endomorphisms of $\sY$ which are
    universal homeomorphisms, such that $\rho_{\sY}$ is
    compatible with $\rho_{\sX}$,  in the sense of  \cite[Def.~1.3
    (ii)]{Rom05}, 
    \emph{i.e.} $\sigma_{m_1m_2}$ is equal to the composition
    \begin{equation} \label{compatibility}
    \begin{tikzpicture}[xscale=5.5,yscale=-1.3,baseline={([yshift=-.5ex]current bounding box.center)}]
       \node (F) at (0,0)
       {$f\circ\rho_{\sY}(m_2)\circ\rho_{\sY}(m_1)$};
       \node (G) at (1,0) {$\rho_{\sX}(m_2)\circ
       f\circ\rho_{\sY}(m_1)$};  
       \node (H) at (2,0)  {$\rho_{\sX}(m_2)\circ\rho_{\sX}(m_1)\circ
    f$};
       \node (I) at (0,1) {$f\circ\rho_{\sY}(m_1m_2)$};
       \node (J) at (2,1) {$\rho_{\sX}(m_1m_2)\circ f$};
       \draw[->] (F) -- (G)
       node[midway,above,scale=.6]{$\sigma_{m_2}$};
       \draw[->] (G) -- (H)
       node[midway,above,scale=.6]{$\sigma_{m_1}$}; 
       \draw[<-] (F) -- (I)
       node[midway,left,scale=.6]{$\cong$};
       \draw[%double equal sign distance,
       ->] (H) -- (J)node[midway,left,scale=.6]{$\cong$};
   \end{tikzpicture}
   \end{equation}
   and $\sigma_{1}$ equals the composition $f\circ\rho_{\sY}(1)\cong
   f\circ\id_{\sY}=f$. The 1-morphisms 
   $(\sY',\rho_{\sY}',\sigma')\to(\sY,\rho_{\sY},\sigma)\in\FEt(\sX/M)$ are
   1-morphisms $g\colon\sY'\to\sY\in \FEt(\sX)$ together
   with 2-morphisms $\delta_m\colon
   \rho_{\sY}(m)\circ g\cong g\circ\rho_{\sY'}(m)$ making $g$ an
   $M$-equivariant map (cf.~\cite[Def.~1.3 (ii)]{Rom05}) and are compatible
   with $\sigma,\sigma'$. The 2-morphisms between two morphisms
      \[\begin{tikzpicture}[xscale=4.17,yscale=2.2,baseline={([yshift=-.5ex]current bounding box.center)}]
             
             \node (A0_1) at (1, 0) {$(\sY',\rho_{\sY}',\sigma')$};
             \node (A0_2) at (2, 0) {$(\sY,\rho_{\sY},\sigma)$};
       
                          \draw[>=latex,->] ([yshift= 1pt] A0_1.east) to
                              [out=10,in=170]
                              node[above,scale=0.5]{$(g,\delta)$}([yshift= 1pt] A0_2.west);
             \draw[>=latex,->] ([yshift=-1pt] A0_1.east) to [out=-10,in=-170]
                 node[below,scale=0.5]{$(g',\delta')$} ([yshift=-1pt] A0_2.west);
             \end{tikzpicture}\] are 2-morphisms in $g\Rightarrow g'\in\FEt(\sX)$ which are
      compatible with $\delta_m,\delta_m'$ for all $m$.
  \end{defn}
      Note that since $\FEt(\sX)$ is anti-equivalent to the category of finite
étale $\sO_{\sX}$-algebras, it is essentially a $1$-category. That is, for any
two
$1$-morphisms $g, g'\colon A\to B$,  the set of 2-morphisms
$\Hom(g,g')$ contains at most one element. Consequently, $\FEt(\sX/M)$ is also
essentially a $1$-category.

However, the introduction of a monoid action in $\FEt(\sX/M)$ — where objects are
finite étale $M$-equivariant morphisms — introduces additional
complexity. Although the resulting category remains a $1$-category, the
equivariance condition inherently involves $2$-morphisms through higher
compatibility requirements. Specifically, one must ensure that these
compatibilities hold at all levels. It turns out, however, that these higher
compatibilities are automatic.

    \begin{lem}\label{universal homeomorphism vs cartesian} For
    any $m\in M$ and $\sY\in\FEt(\sX)$ set
    $m^*\sY\coloneqq \sY\times_{\rho_{\sX}(m)}\sX$. The following
        categories are equivalent
        \begin{enumerate}[label={\rm (\alph*)}]
            \item the
        category
        $\FEt(\sX/M)$;
            \item  the category of  2-cartesian diagrams
        \begin{equation}\label{cartesian blocks}
            \begin{tikzpicture}[xscale=1.5,yscale=0.9,bmr/.pic={\draw (0,0)--++(-90:2mm)--++(180:2mm);},baseline={([yshift=-.5ex]current bounding box.center)}]
                 \path
                        (0,0)     node (F) {$\sY$}
                     +(-43:.4) pic[scale=1,red]{bmr}
                     +(0:1.5)  node (star) {$\sY$}
                     ++(-90:1.5) node (X) {$\sX$}
                     +(0:1.5)  node (Y) {$\sX$};
                     \node (start) at (1.1,-0.5) {};
                     \node (end) at (.5,-1.1) {};
                 \draw[->] (F)--(star) node[midway,above,scale=.6]{$\rho_{\sY}(m)$};
                 \draw[->] (F)--(X) node[midway,left,scale=.6]{$f$};
                 \draw[->] (X)--(Y) node[midway,above,scale=.6]{$\rho_{\sX}(m)$};
                 \draw[->] (star)--(Y) node[midway,right,scale=.6]{$f$};
                 \draw[double distance=3.5pt,line width=.8pt,-{Classical TikZ
                 Rightarrow[length=3mm, line width=.8pt]}] (start)--(end)
                 node[midway,above,scale=.6,yshift=7pt]{$\sigma_m$};
            \end{tikzpicture}
        \end{equation}
        where $\sY\xrightarrow{f}\sX\in\FEt(\sX)$ and $\rho_{\sY}$ is an
        action of $M$ on $\sY$ via
        $\F_q$-endomorphisms 
        %(not necessarily universal homeomorphisms) 
        and $\sigma_m$ is the 2-isomorphism
        as in the definition of $\FEt(\sX/M)$;
    \item the category of tuples $(\sY,\alpha_{\sY},\iota,\tau)$,
        where $\sY\in\FEt(\sX)$, $\alpha_{\sY,m}\colon
        \sY\xrightarrow{\simeq} m^*\sY\in\FEt(\sX)$ is an
        equivalence for each $m\in M$, $\iota_{m_1,m_2}\colon
        m_1^*\alpha_{\sY,m_2}\circ\alpha_{\sY,m_1}\cong
        \alpha_{\sY,m_1m_2}$ and $\tau\colon \alpha_{\sY,1}\cong \id_{\sY}$
        are 2-isomorphisms in $\FEt(\sX)$. They satisfy the usual
        associativity and unit diagrams for a monoid action. A morphism
        $g\colon\sY'\to\sY$ is equipped with 2-isomorphisms
        $\theta_m\colon\alpha_{\sY,m}\circ g\cong
        m^*g\circ\alpha_{\sY',m}$ compatible with $\iota$ and $\tau$;
    \item the category of pairs $(\sY,\alpha_{\sY})$, where
        $\sY\in\FEt(\sX)$, $\alpha_{\sY,m}\colon
        \sY\xrightarrow{\simeq} m^*\sY\in\FEt(\sX)$ is an
        equivalence for each $m\in M$. Moreover, we have $m_1^*\alpha_{\sY,m_2}\circ\alpha_{\sY,m_1}\cong
        \alpha_{\sY,m_1m_2}$ and $\alpha_{\sY,1}\cong \id_{\sY}$. A morphism
        $(\sY',\alpha_{\sY'})\to(\sY,\alpha_{\sY})$ is a map $g\colon
        \sY'\to\sY\in\FEt(\sX)$ commuting with $\alpha_{\sY'}$ and $\alpha_{\sY}$. 
    \item the category of pairs $(\sA, \alpha_{\sA})$, where
        $\sA$ is a finite étale $\sO_{\sX}$-algebra and $\alpha_{\sA,m}\colon
        m^*\sA\to \sA$ $(m^*\sA\coloneqq \rho_{\sX}(m)^*\sA)$ is an
        $\sO_{\sX}$-algebra isomorphism satisfying
        $\alpha_{\sA,e}=\id_{\sA}$ and
        $\alpha_{\sA,m_1m_2}=\alpha_{\sA,m_1}\circ m_1^*\alpha_{\sA,m_2}$. A morphism
        $(\sA',\alpha_{\sA'})\to(\sA,\alpha_{\sA})$ is a map $g\colon
        \sA'\to\sA$ of $\sO_{\sX}$-algebras commuting with $\alpha_{\sA'}$ and
        $\alpha_{\sA}$.
        \end{enumerate}
\end{lem}
\begin{proof} First of all, the categories in (a), (b), (c) and (d) are essentially
    1-categories as remarked before. Thus we can ignore the natural 2-morphisms
    in them.
    Let's first show (a)$\Leftrightarrow$(b).  
    Consider the 2-commutative diagram:
    \begin{equation}\label{M-action via automorphisms and endomorphisms}
        %%%%%% Arrow tip style [-stealth], [-to], [-latex], [-|], [>=latex,->], [|->]
        %%%%%% Arrow style: [red], [dashed], [dotted, line width=12pt], [ultra thick], [thick], [semithick], [thin], [very thin], [ultra thin]
        %%%%%% Double arrows: \draw ([yshift=-2pt] A0_0.east) -- ([yshift=-2pt] A0_1.west);
        %%%%%% Edge follow the rules of draw.
        %%%%%% Node style: [above,scale=0.9], [below], [midway], [anchor=south], [ellipse,draw], [circle,fill=red!20] 
        %%%%%% Arrow Node: [pos=0.5] (this is the default position)
        %%%%%% Equation tag position center alignment: [baseline=(current  bounding  box.center)]
            \begin{tikzpicture}[xscale=2.0,yscale=1.2,bmr/.pic={\draw
                (0,0)--++(-90:2mm)--++(180:2mm);},baseline={([yshift=-.5ex]current bounding box.center)}]
                 \path
                        (0,0)     node (F) {$m^*\sY$}
                     +(-43:.4) pic[scale=1,red]{bmr}
                     +(0:1.5)  node (star) {$\sY$}
                     ++(-90:1.5) node (X) {$\sX$}
                     +(0:1.5)  node (Y) {$\sX$}
                     ++(90:2.5)
                     +(-180:1.) node (Z) {$\sY$};
                 \node (start) at (1.3,-0.25) {};
                 \node (end) at (.2,-1.35) {};
                 \draw[->] (F)--(star);
                 \draw[->] (F)--(X);
                 \draw[->] (X)--(Y) node[midway,above,scale=.6]{$\rho_{\sX}(m)$};
                 \draw[->] (star)--(Y);
                 \draw[dashed,->] (Z)--(F)
                 node[midway,above,scale=.6]{$\hspace{10pt}\alpha_{\sY,m}$};
                 \draw[dotted,->] (Z) to [out=30,in=140] node[midway,above,scale=.6]{$\hspace{10pt}\rho_{\sY}(m)$} (star);
                 \draw[->] (Z) to [out=-90,in=140](X);
                 \draw[double distance=3.5pt,line width=.8pt,-{Classical TikZ
                 Rightarrow[length=3mm, line width=.8pt]}] (start)--(end);
            \end{tikzpicture}
    \end{equation}
    Suppose the square defined by $\rho_{\sY}(m)$ and
    $\rho_{\sX}(m)$ is cartesian, then $\alpha_{\sY,m}$ is an
    equivalence, so $\rho_{\sY}(m)$, as a pullback of $\rho_{\sX}(m)$, is a
    universal homeomorphism, \emph{i.e.} the square belongs to  $\FEt(\sX/M)$.
    Conversely, if
        $\rho_{\sY}(m)$ is a universal homeomorphism, then $\alpha_{\sY,m}$ is
        an isomorphism. Indeed, if $\rho_{\sY}(m)$ is
        a universal homeomorphism, then the fact that $m^*\sY\to\sY$ is a
        universal homeomorphism implies that $\alpha_{\sY,m}$ is a universal
        homeomorphism on the nose (cf.~Lemma \ref{decent geometric properties
        of uni homeomorphisms}).  Since both $\sY$ and $m^*\sY$ are in
        $\FEt(\sX)$, $\alpha_{\sY,m}$ is an isomorphism at
        all geometric fibers of $\sX$, so it is a degree 1 finite étale
        map, \emph{i.e.} an equivalence, hence the square defined by
        $\rho_{\sY}(m)$ and $\rho_{\sX}(m)$ is
        cartesian. 

        (b)$\Leftrightarrow$(c) If we ignore the action, then the data of diagram
        \eqref{cartesian blocks} are equivalent to  the data of  pairs
        $(\sY,\alpha_{\sY})$, where $\alpha_{\sY,m}\colon \sY\to m^*\sY$ is a
        1-morphism for each $m\in M$. The data of the 2-morphisms $\iota,\tau$ correspond
        to the data of 2-morphisms defining the action $\rho_{\sY}$
        (cf.~\cite[Def. 1.3 (i)]{Rom05}) and the compatibility condition
        defined in diagram \eqref{compatibility} and the line after. The higher
        compatibilities of $\iota,\tau$ correspond to the higher associativity
        and other higher compatibility constraints on the action
        $\rho_{\sY}$.

(c)$\Leftrightarrow$(d) For this we just have to note that the higher
compatibilities are automatic because $\iota,\tau$ are 2-morphisms in
$\FEt(\sX)$ which is essentially a 1-category.

(d)$\Leftrightarrow$(e)
        This follows immediately from the
        anti-equivalence $\FEt(\sX)\simeq$ \{finite étale
        $\sO_{\sX}$-algebras\}. 
    \end{proof}

\begin{rmk} From Lemma \ref{universal homeomorphism vs cartesian} (e) one sees that $\FEt(\sX/M)$ is just the category $\FEt(\sX)$ rigidified along a compatible $M$-collection of maps in $\FEt(\sX)$. It follows that $\FEt(\sX/M)$ ``inherits'' important categorical features of $\FEt(\sX)$, e.g. the functor \resizebox{4cm}{!}{$\FEt(\sX/M)\to \FEt(\sY/M)$} induced by a map  \resizebox{3.3cm}{!}{$f\colon(\sY,\rho_{\sY})\to(\sX,\rho_{\sX})$} is fully faithful (resp. equivalence), if \resizebox{3cm}{!}{$\FEt(\sX)\to \FEt(\sY)$} is so. Moreover,  $f$ is a map of (effective) descent for $\FEt(-/M)$ if it is so for $\FEt(-)$.
\end{rmk}

An important consequence of Lemma \ref{universal homeomorphism vs cartesian} is
the following main theorem of this section. 
\begin{thm} \label{Galois category} If $\sX$ is $G$-connected, then
    $\FEt(\sX/M)$ is a Galois category.
\end{thm}
\begin{proof}
With the above preparation, one could verify the axioms listed in
\cite[\href{https://stacks.math.columbia.edu/tag/0BMY}{0BMY}]{stacks-project}
sequentially, following the same method as in
\cite[\href{https://stacks.math.columbia.edu/tag/0BL6}{0BL6}]{stacks-project}.
Alternatively, one could verify the axioms listed in \cite[Exposé V, \S4]{SGA1}
sequentially, following the same method as  in \cite[Thm. 4.2]{Noohi04}.
\end{proof}
To ensure future flexibility, we provide an alternative proof of Theorem \ref{Galois category} within the topos-theoretic framework of \cite{Leroy79,Zoo01}, as suggested by an anonymous referee.  Nonetheless, we shall warn the
reader that the proof only works under Setting \ref{the assumption for the
topos approach} below, which includes the cases of primary interest to us.

\begin{defn}
    The \textit{Weil-lisse-étale site}, denoted $(\sX,\rho_{\sX})_{\fppf}$, consists of:
\begin{itemize}
\item \textbf{Objects}: Smooth $M$-equivariant morphisms $\alpha_U\colon
    (U,\rho_U) \to (\sX,\rho_{\sX})$, where:
\begin{itemize}
\item $U$ is an $\F_q$-algebraic space,
\item The $M$-action $\rho_U$ on $U$ is via universal homeomorphisms
\end{itemize}
\item \textbf{Morphisms}: Pairs $(f, \sigma_{U_1,U_2})$ where:
\begin{itemize}
\item $f\colon (U_1,\rho_{U_1}) \to (U_2,\rho_{U_2})$ is an $M$-equivariant morphism,
\item $\sigma_{U_1,U_2}\colon \alpha_{U_2}\circ f \cong \alpha_{U_1}$ is a 2-isomorphism
\end{itemize}
\item \textbf{Coverings}: Families of étale $M$-equivariant morphisms that are jointly surjective
\end{itemize}
\end{defn}

At this stage, it is not yet clear whether the category $(\sX,\rho_{\sX})_{\fppf}$ 
contains any nonempty objects. To ensure that our definition is useful, we introduce 
the following assumption:

\begin{setting}\label{the assumption for the topos approach}
    There exists an atlas $u\colon X \to \sX$ that lifts the $M$-action. 
    More precisely, we assume the existence of a representable, smooth, and surjective 
    morphism $u\colon (X,\rho_{X})\to(\sX,\rho_{\sX})$ in the category 
    $(\sX,\rho_{\sX})_{\fppf}$.
\end{setting}

This condition holds in the following key cases:
\begin{itemize}
    \item {\bf Product of stacks with partial Frobenius action:} The setting of \S 6, where $\sX = \sX_1 \times \cdots \times \sX_n$ and $\rho_{\sX}$ is given by partial Frobenii.

    \item {\bf Algebraic spaces:} When $\sX = X$ is itself an algebraic space.

    \item {\bf Trivial action:} The case where $M$ is trivial.
\end{itemize}

With this assumption in mind, we can commence by noting several defining features
of this Weil-lisse-étale site.
    
\begin{lem}\label{subcanonical site}
The category $(\sX,\rho_{\sX})_{{\fppf}}$, equipped with the specified coverings, forms a subcanonical site. That is, all representable presheaves on this site are sheaves.
\end{lem}

The verification of the site axioms \cite[\href{https://stacks.math.columbia.edu/tag/00VH}{00VH}]{stacks-project}, as well as the subcanonicity condition, is straightforward; the presence of the $M$-action introduces no additional difficulty.

Let's consider the topos $\sT\coloneqq \Sh((\sX,\rho_{\sX})_{\fppf})$. Recall that a
            sheaf on a site is called \textit{finite locally constant} if it is locally isomorphic to a constant sheaf associated to a finite set
            (cf.~\cite[\href{https://stacks.math.columbia.edu/tag/093P}{093P}]{stacks-project}, \cite[Chap. I, Déf. 1.1]{Zoo01}). 
           For any object $(\sY,\rho_{\sY})\in\FEt(\sX/M)$, we define a presheaf $h_{(\sY,\rho_{\sY})}$
           ``represented'' by this object via the assignment:
           \[ 
               (U,\rho_U)\longmapsto \Hom_{(\sX,\rho_{\sX})}((U,\rho_U),(\sY,\rho_{\sY}))
           \] 
           Note that for any object $(U,\rho_{U})$ of
           $(\sX,\rho_{\sX})_{\fppf}$,
           $h_{(\sY,\rho_{\sY})}|_{(U,\rho_{U})}=h_{(\sY_U,\rho_{\sY_U})}$ is a sheaf by Lemma \ref{subcanonical site}. It follows that
           $h_{(\sY,\rho_{\sY})}$ is a sheaf, and we thereby obtain a
           well-defined functor \[\FEt(\sX/M)\longrightarrow\sT\]

\begin{lem}               
    The functor $\FEt(\sX/M) \to \sT$ is fully faithful, and its essential image 
    consists exactly of the finite locally constant sheaves.
\end{lem}

\begin{proof}
    We first show that the essential image of $\FEt(\sX/M)$ is contained in 
    $\sT_{\flc}$, the full subcategory of $\sT$ consisting of finite locally constant sheaves.
    
    Let $f\colon (\sY,\rho_{\sY})\to(\sX,\rho_{\sX})$ be finite \'{e}tale.
    Its rank loci are open and closed and are preserved by the $M$-action.
    Pulling them back to the equivariant atlas in Setting
    \ref{the assumption for the topos approach} gives a jointly surjective
    family of objects of the site. We may therefore work on one rank locus
    and assume that $f$ has constant degree $n$; no boundedness of the degree
    on the original, possibly non-quasi-compact stack is required.
    The diagonal map 
    $(\sY,\rho_{\sY})\xrightarrow{\Delta}(\sY,\rho_{\sY})\times_{(\sX,\rho_{\sX})}(\sY,\rho_{\sY})$ 
    is a clopen immersion, so its complement
    $(\sY,\rho_{\sY})\times_{(\sX,\rho_{\sX})}(\sY,\rho_{\sY})\setminus\Delta(\sY,\rho_{\sY})$ 
    is a finite étale morphism of degree $n-1$. 
    By induction, this complement splits over some lisse-étale covering 
    $(\sZ,\rho_{\sZ})\to(\sY,\rho_{\sY})$; moreover, this same covering 
    $(\sZ,\rho_{\sZ})\to(\sX,\rho_{\sX})$ also splits $(\sY,\rho_{\sY})$.
    
    To complete the argument, it suffices to note that both $\FEt(\sX/M)$ and $\sT_{\flc}$ 
    satisfy Weil-lisse-étale descent, and they agree on constant objects. 
    For $\FEt(\sX/M)$, descent follows from Lemma \ref{universal homeomorphism vs cartesian} (e) 
    together with flat descent for quasi-coherent sheaves. 
    For $\sT_{\flc}$, this is a general property of the topos (see \cite[Example 4.11]{Vis05}).
\end{proof} 
         
     By \cite[\S 5.2.2]{Leroy79} or \cite[\S 8.40-8.50]{Johnstone77}, the category of finite locally constant sheaves in a connected topos forms a Galois category. Therefore, to establish Theorem \ref{Galois category}, it suffices to verify that the topos $\Sh((\sX,\rho_{\sX})_{\fppf})$ is connected when $\sX$ is $G$-connected.

\begin{proof}[Proof of Theorem \ref{Galois category} under Setting \ref{the assumption for the topos approach}]
    We show that $\Sh((\sX,\rho_{\sX})_{\fppf})$ is a connected topos by proving that its terminal object $*$ is connected.
    
    Suppose, to the contrary, that $*$ decomposes as a coproduct $* = \sF_1 \sqcup \sF_2$ with $\sF_1, \sF_2$ nonempty. The maps
    $\sF_i\to *$ are complementary subterminal objects. Although coproducts
    of sheaves involve sheafification, the epimorphism
    $\sF_1\sqcup\sF_2\to *$ means precisely that its sections lift locally.
    For $i = 1,2$, define
    \[
    S_i \coloneqq \left\{ f\colon (U,\rho_U)\to(\sX,\rho_{\sX}) \; \middle| \;
    f \in (\sX,\rho_{\sX})_{\fppf},\, U\neq\emptyset, \text{ and } \sF_i(f) \neq \emptyset \right\}.
    \]
    
    We claim that if $f \in S_1$ and $g \in S_2$, then $f(U) \cap g(V) = \emptyset$.
    Moreover, the two unions below cover $|\sX|$: apply the preceding local
    lifting property after pulling back to an equivariant atlas supplied by
    Setting \ref{the assumption for the topos approach}. Since smooth maps
    are open and all maps in the site are equivariant, we obtain a
    $G$-equivariant decomposition into two nonempty open subsets
    \[
        |\sX| = \bigcup_{f \in S_1} f(U) \ \sqcup \ \bigcup_{g \in S_2} g(V),
    \]
    which contradicts the assumption that $\sX$ is $G$-connected.

    To prove the claim, suppose that $f(U) \cap g(V) \neq \emptyset$. Then the fiber product $W \coloneqq U \times_{\sX} V$ is nonempty.
    Equip $W$ with the $M$-action $\rho_W \coloneqq \rho_U \times_{\rho_{\sX}} \rho_V$.

    We first verify that $(W,\rho_W)$ belongs to $(\sX,\rho_{\sX})_{\fppf}$.
    For each $m \in M$, the morphism $\rho_W(m)$ is a universal homeomorphism on the nose by Lemma \ref{decent geometric properties of uni homeomorphisms}.
    Since $\rho_W(m)$ is a morphism between algebraic spaces, Lemma \ref{decent geometric properties of uni homeomorphisms} (3) implies that it is actually a universal homeomorphism.

    The unique map $h_{(U,\rho_U)} \to *$ factors through $\sF_1$, and similarly $h_{(V,\rho_V)} \to *$ factors through $\sF_2$.
    Consequently, their fiber product $h_{(W,\rho_W)} \cong h_{(U,\rho_U)} \times_{*} h_{(V,\rho_V)}$ must factor through $\sF_1 \times_* \sF_2$.
    But $h_{(W,\rho_W)}$ is nonempty (since $W$ is nonempty), while
    $\sF_1\times_*\sF_2$ is initial because the two subobjects are
    complementary.
    This is a contradiction, which establishes the claim.
\end{proof}

\section{Partial Frobenius maps for étale covers}\label{Drinfeld's lemma} Let $\sX_1,\sX_2,\dots,\sX_n$ be algebraic stacks over $\F_q$, and set $\sX\coloneqq 
\sX_1\times_{\F_q}\sX_2\times_{\F_q}\cdots\times_{\F_q}\sX_n$. Denote
$\phi_{i}\coloneqq \phi_{\sX_i}$
the $q$-Frobenius of $\sX_i$ for $1\leq i\leq n$. Then $\sX$ is
equipped with a strict
$M\coloneqq\phi_1^\N\times\phi_2^\N\times\cdots\times\phi_n^\N$-action. Let $M_0$ denote the submonoid
 $\phi_1^\N\times\cdots\times\phi_{n-1}^\N$, and let $G$ (resp. $G_0$)
denote
the free abelian group generated by $M$ (resp. $M_0$). The category
$\FEt(\sX/M_0)$ can be understood as 2-cartesian diagrams
\begin{equation}\label{n-1}
           \begin{tikzpicture}[xscale=1.5,yscale=0.9,bmr/.pic={\draw (0,0)--++(-90:2mm)--++(180:2mm);},baseline={([yshift=-.5ex]current bounding box.center)}]
                 \path
                        (0,0)     node (F) {$\sY$}
                     +(-43:.4) pic[scale=1,red]{bmr}
                     +(0:1.5)  node (star) {$\sY$}
                     ++(-90:1.5) node (X) {$\sX$}
                     +(0:1.5)  node (Y) {$\sX$};
                 \path
                        (star)
                     +(-43:.4) pic[scale=1,red]{bmr}
                     +(0:1.5)  node (star1) {$\sY$}
                     ++(-90:1.5) 
                     +(0:1.5)  node (Y1) {$\sX$};
                 \path
                        (star1)
                     +(-43:.4) pic[scale=1,red]{bmr}
                     +(0:1.5)  node (star2) {$\sY$}
                     ++(-90:1.5) 
                     +(0:1.5)  node (Y2) {$\sX$};
                 \path
                        (star2)
                     +(-43:.4) pic[scale=1,red]{bmr}
                     +(0:1.5)  node (star3) {$\sY$}
                     ++(-90:1.5) 
                     +(0:1.5)  node (Y3) {$\sX$};

                     \node (start) at (1.1,-0.5) {};
                     \node (end) at (.5,-1.1) {};
                     \node (start1) at (2.6,-0.5) {};
                     \node (end1) at (2,-1.1) {};
                     \node (start2) at (4.1,-0.5) {};
                     \node (end2) at (3.5,-1.1) {};
                     \node (start3) at (5.6,-0.5) {};
                     \node (end3) at (5,-1.1) {};

                 \draw[->] (F)--(star)
                 node[midway,above,scale=.6]{$\varphi_{n-1}$};
                 \draw[->] (F)--(X) node[midway,left,scale=.6]{$f$};
                 \draw[->] (X)--(Y)
                 node[midway,above,scale=.6]{$\phi_{n-1}$};
                 \draw[->] (star)--(Y) node[midway,left,scale=.6]{$f$};
                 \draw[->] (star)--(star1)
                 node[midway,above,scale=.6]{$\varphi_{n-2}$};
                 \draw[->] (star1)--(Y1) node[midway,left,scale=.6]{$f$};
                 \draw[->] (Y)--(Y1)
                 node[midway,above,scale=.6]{$\phi_{n-2}$};

                 \draw[->] (star1)--(star2)
                 node[midway,above,scale=.6]{$\cdots$};
                 \draw[->] (star2)--(Y2) node[midway,left,scale=.6]{$f$};
                 \draw[->] (Y1)--(Y2)
                 node[midway,above,scale=.6]{$\cdots$};

                 \draw[->] (star2)--(star3)
                 node[midway,above,scale=.6]{$\varphi_{1}$};
                 \draw[->] (star3)--(Y3) node[midway,right,scale=.6]{$f$};
                 \draw[->] (Y2)--(Y3)
                 node[midway,above,scale=.6]{$\phi_{1}$};

                 \draw[double distance=3.5pt,line width=.8pt,-{Classical TikZ
                 Rightarrow[length=3mm, line width=.8pt]}]
                 (start1)--(end1)
                 node[midway,above,scale=.6,yshift=7pt]{$\sigma_{n-2}$};   
                 \draw[double distance=3.5pt,line width=.8pt,-{Classical TikZ
                 Rightarrow[length=3mm, line width=.8pt]}] (start)--(end)
                 node[midway,above,scale=.6,yshift=7pt]{$\sigma_{n-1}$};
                 \draw[double distance=3.5pt,line width=.8pt,-{Classical TikZ
                 Rightarrow[length=3mm, line width=.8pt]}]
                 (start3)--(end3)
                 node[midway,above,scale=.6,yshift=7pt]{$\sigma_1$};
                 \draw[dotted] (start2) -- (end2);
            \end{tikzpicture}
\end{equation}where the $\varphi_i$s (together with the $\sigma_i$s) mutually commute. Morally,
$\FEt(\sX/M_0)$ is just a collection of commuting square blocks. 

\begin{lem}\label{are mutually commute} The category $\FEt(\sX/M_0)$ can be reinterpreted as 2-cartesian
    diagrams
    \begin{equation}\label{n}
           \begin{tikzpicture}[xscale=1.5,yscale=0.9,bmr/.pic={\draw (0,0)--++(-90:2mm)--++(180:2mm);},baseline={([yshift=-.5ex]current bounding box.center)}]
                 \path
                        (0,0)     node (F) {$\sY$}
                     +(-43:.4) pic[scale=1,red]{bmr}
                     +(0:1.5)  node (star) {$\sY$}
                     ++(-90:1.5) node (X) {$\sX$}
                     +(0:1.5)  node (Y) {$\sX$};
                 \path
                        (star)
                     +(-43:.4) pic[scale=1,red]{bmr}
                     +(0:1.5)  node (star1) {$\sY$}
                     ++(-90:1.5) 
                     +(0:1.5)  node (Y1) {$\sX$};
                 \path
                        (star1)
                     +(-43:.4) pic[scale=1,red]{bmr}
                     +(0:1.5)  node (star2) {$\sY$}
                     ++(-90:1.5) 
                     +(0:1.5)  node (Y2) {$\sX$};
                 \path
                        (star2)
                     +(-43:.4) pic[scale=1,red]{bmr}
                     +(0:1.5)  node (star3) {$\sY$}
                     ++(-90:1.5) 
                     +(0:1.5)  node (Y3) {$\sX$};

                     \node (start) at (1.1,-0.5) {};
                     \node (end) at (.5,-1.1) {};
                     \node (start1) at (2.6,-0.5) {};
                     \node (end1) at (2,-1.1) {};
                     \node (start2) at (4.1,-0.5) {};
                     \node (end2) at (3.5,-1.1) {};
                     \node (start3) at (5.6,-0.5) {};
                     \node (end3) at (5,-1.1) {};

                 \draw[->] (F)--(star)
                 node[midway,above,scale=.6]{$\varphi_{n}$};
                 \draw[->] (F)--(X) node[midway,left,scale=.6]{$f$};
                 \draw[->] (X)--(Y)
                 node[midway,above,scale=.6]{$\phi_{n}$};
                 \draw[->] (star)--(Y) node[midway,left,scale=.6]{$f$};
                 \draw[->] (star)--(star1)
                 node[midway,above,scale=.6]{$\varphi_{n-1}$};
                 \draw[->] (star1)--(Y1) node[midway,left,scale=.6]{$f$};
                 \draw[->] (Y)--(Y1)
                 node[midway,above,scale=.6]{$\phi_{n-1}$};

                 \draw[->] (star1)--(star2)
                 node[midway,above,scale=.6]{$\cdots$};
                 \draw[->] (star2)--(Y2) node[midway,left,scale=.6]{$f$};
                 \draw[->] (Y1)--(Y2)
                 node[midway,above,scale=.6]{$\cdots$};

                 \draw[->] (star2)--(star3)
                 node[midway,above,scale=.6]{$\varphi_{1}$};
                 \draw[->] (star3)--(Y3) node[midway,right,scale=.6]{$f$};
                 \draw[->] (Y2)--(Y3)
                 node[midway,above,scale=.6]{$\phi_{1}$};

                 \draw[double distance=3.5pt,line width=.8pt,-{Classical TikZ
                 Rightarrow[length=3mm, line width=.8pt]}]
                 (start1)--(end1)
                 node[midway,above,scale=.6,yshift=7pt]{$\sigma_{n-1}$};   
                 \draw[double distance=3.5pt,line width=.8pt,-{Classical TikZ
                 Rightarrow[length=3mm, line width=.8pt]}] (start)--(end)
                 node[midway,above,scale=.6,yshift=7pt]{$\sigma_{n}$};
                 \draw[double distance=3.5pt,line width=.8pt,-{Classical TikZ
                 Rightarrow[length=3mm, line width=.8pt]}]
                 (start3)--(end3)
                 node[midway,above,scale=.6,yshift=7pt]{$\sigma_1$};
                 \draw[dotted] (start2) -- (end2);
            \end{tikzpicture}
        \end{equation}where the $\varphi_i$s {\rm (}together with the
        $\sigma_i$s{\rm )} mutually commute and
        $\varphi_1\circ\varphi_2\circ\cdots\circ\varphi_n\cong \phi_{\sY}$.
\end{lem}
\begin{proof} Note that given diagram \eqref{n-1}, there is a unique way to
    complete it to diagram \eqref{n} thanks to the universality of 2-cartesian
    diagrams. One still has to show that $\varphi_n$ commutes with $\varphi_i$
    for $1\leq i\leq n-1$. By induction and the universality again, we only
    have to show it for $n=2$. Suppose $\varphi_1\circ\varphi_2\cong \phi_{\sY}$, then
    we have
    \begin{equation}\label{natural isomorphisms}
           \varphi_1\circ\varphi_2\circ\varphi_1\cong
           \phi_{\sY}\circ\varphi_1\cong \varphi_1\circ \phi_{\sY}
    \end{equation}
    Now the universality of 2-cartesian diagrams provides the ``left
    cancellation'', thus $\varphi_2\circ\varphi_1\cong \phi_{\sY}$. As
    all the isomorphisms in \eqref{natural isomorphisms} are natural, the
    isomorphism
    $\varphi_2\circ\varphi_1\cong \phi_{\sY}\cong \varphi_1\circ\varphi_2$ is naturally compatible with
    $\sigma_1,\sigma_2$. 
\end{proof}

\begin{rmk}\label{dropping out any i}
Note that we have dropped $\varphi_n$ in the definition of $\FEt(\sX/M_0)$, but
in view of Lemma \ref{are mutually commute}
dropping
any other $\phi_i$  would have yielded the same category. 
\end{rmk}

\section{The permanence of $\pi_1^\et$}

As discussed in the introduction, the permanence of $\pi_1^\et$ played a
pivotal role in the proof of Drinfeld's lemma for schemes. Here is
the stack version of it.

\begin{lem}\label{Permanence} Let $\sX$ be any category fibered in groupoid over $\F_q$, and
    let $k$ be any algebraically closed field containing $\F_q$. Set $\sX_k\coloneqq
    \sX\times_{\F_q}k$ and $\phi_k\colon \sX_k\to\sX_k$ the partial Frobenius
    of the second factor. Then the pullback along $\sX_k\to\sX$ induces an equivalence
    \[
        \FEt(\sX)\xrightarrow{\quad \simeq\quad}\FEt\left(\sX_k/\phi_k^\N\right)
    \]
\end{lem}
\begin{proof} The simplest approach is to apply \cite[Theorem II (C)]{DTZ20}, observing that the pullback functor induces an equivalence of categories between finite étale covers on both sides.

    One could
    alternatively apply the simple independent lemma \cite[Lemma 3.2]{DTZ20}  to reduce the problem to the case
    where $\sX$ is an affine scheme, then conclude with \cite[Lemma 4.2.6]{Ked17} or \cite[Lemma 16.2.6]{SW20}.
\end{proof}

Let $\sX$ be an algebraic stack over $\F_q$ equipped with a monoid action $\rho_{\sX}\colon M \to \End_{\F_q}(\sX)$ via universal homeomorphisms. Let $G$ denote the Grothendieck group associated to $M$. We now state a slightly generalized version of Lemma \ref{Permanence} that will be essential for our subsequent arguments.

\begin{cor}\label{Permanence v2} The pullback along $\sX_k\to\sX$ induces an equivalence:
    \begin{equation}\label{Drinfeld-Lau descent for G-covers}
\FEt(\sX/M)\xrightarrow{\quad
\simeq\quad}\FEt\left(\sX_k/(M\times\phi_k^\N)\right)
\end{equation}
\end{cor}
\begin{proof} In view of Lemma \ref{Permanence}, there is an equivalence $\FEt(\sX)\xrightarrow{
\simeq}\FEt\left(\sX_k/\phi_k^\N\right)
$. An object of the left hand side of \eqref{Drinfeld-Lau descent for G-covers} can be
interpreted, via Lemma \ref{universal homeomorphism vs cartesian} (d), as a
collection of morphisms $\left\{\sigma_m\colon\sY\to m^*\sY\right\}_{m\in
M}$
in $\FEt(\sX)$  which are
compatible in a natural way. An object of the right hand side of \eqref{Drinfeld-Lau descent for
G-covers} can be interpreted similarly as a collection
$\left\{\delta_m\colon\sZ\to m^*\sZ\right\}_{m\in
M}$
 in $\FEt(\sX_k)$
 plus an action of $\phi_k^{\N}$ on $\sZ$: $\{f_n\colon \sZ\to
 \phi_k^{n*}\sZ\}_{n\in\N}$. Since the $M$-action commutes with the
$\phi_k^\N$-action, we have the following 2-commutative diagram

\[
\xymatrix{
    \sZ \ar[rr]^{f_n} \ar[d]_{\delta_m} && \phi_k^{n*}\sZ
    \ar[d]^{\phi_k^{n*}\delta_m} \\
    m^*\sZ \ar[rr]^{m^*f_n} && \phi_k^{n*}m^*\sZ
}
\]
for $\forall m\in M$, $\forall n\in\N$. This shows that the morphisms
$\delta_m$ are $\phi_k^\N$-equivariant.  The assertion then follows from the equivalence $\FEt(\sX)\xrightarrow{
\simeq}\FEt\left(\sX_k/\phi_k^\N\right)
$.
\end{proof}

\section{Drinfeld's lemma}
We are now ready to prove Theorem \ref{Drinfeld's Lemma}. As a key ingredient,
we first establish a Drinfeld-style analogue of Grothendieck's existence
theorem, where the usual properness assumption is replaced by partial Frobenius
actions. We adopt the notation introduced in \S6.

\textbf{Setting.}
Throughout this section, 
$\sX_1,\sX_2,\dots,\sX_n$ are algebraic stacks over $\F_q$, and we set
$\sX = \sX_1\times_{\F_q}\cdots\times_{\F_q}\sX_n$.
For each $i$, let $\phi_i$ denote the $q$-Frobenius of $\sX_i$.
Then $\sX$ carries a strict action $\rho_{\sX}$ of the monoid
$M = \phi_1^\N\times\cdots\times\phi_n^\N$.
We also use the following notation:
\(
    G   \coloneqq \phi_1^\Z\times\cdots\times\phi_n^\Z,\,
    G_0 \coloneqq \phi_1^\Z\times\cdots\times\phi_{n-1}^\Z,\,
    M_0 \coloneqq \phi_1^\N\times\cdots\times\phi_{n-1}^\N.
\)
\begin{lem}\label{pi_1 existence theorem}
    Suppose that $\sX_1,\ldots,\sX_n$ are qcqs. Let $Y$ be an $\F_q$-scheme. For any geometric point $\bar{x}\colon\Spec(k)\to Y$ the base change functor induces an equivalence:
    \begin{equation}\label{existence theorem equation}
        \varinjlim_{(U,\bar{u})}
        \FEt(\sX\times_{\F_q}U/M)\longrightarrow
        \FEt(\sX\times_{\F_q}\bar{x}/M)
    \end{equation}
    where $(U,\bar{u})$ runs over all affine étale neighborhoods of $\bar{x}$.
\end{lem}

\begin{proof} 
\textbf{Step 0, general principle.} 
    Let $u\colon\sT'\to\sT$ be a $1$-morphism of stacks over a site
    $\sC$, and set
    $\sT''\coloneqq\sT'\times_{\sT}\sT'$,
    $\sT'''\coloneqq\sT''\times_{\sT}\sT'$. 
    Suppose that $A,B$ are pseudo-functors over $\sC$ (cf.~\cite[Def.
    3.10]{Vis05}) or equivalently fibered categories over $\sC$
    (cf.~\cite[\S 3.1.3]{Vis05}). Denote by $A(\sT)\coloneqq \Hom_{\sC}(\sT,A)$
    the Hom-category, and similarly for the other objects. Suppose
    we have $2$-commutative diagrams 
    \begin{equation}\label{descent on X1}
                    \begin{tikzpicture}[xscale=2.5,yscale=1.2,baseline=(current  bounding  box.center)]
                 \path
                 (0,0)     node[scale=.8] (F)
                 {$A(\sT)$}
                     ++(0:1.5)  node[scale=.8] (star)
                     {$A(\sT')$}
                     ++(0:1.5) node[scale=.8] (ver13)
                     {$A(\sT'')$}
                     ++(0:1.5) node[scale=.8] (ver14)
                     {$A(\sT''')$}
                     ++(-90:1.5) node[scale=.8] (new)
                     {$B(\sT''')$}
 
                     ++(-180:1.5) node[scale=.8] (X)
                     {$B(\sT'')$}
                     ++(-180:1.5)  node[scale=.8] (Y)
                     {$B(\sT')$}
                     ++(-180:1.5) node[scale=.8] (Z)
                     {$B(\sT)$};
                     \draw[->] (F)--node[above,scale=0.6]{$\scriptstyle{u^*}$}(star);
                 \draw[->]
                 ([yshift=-2pt]star.east)--node[below,scale=0.6]{$\scriptstyle{p_2^*}$}([yshift=-2pt]ver13.west);
\draw[->] ([yshift=2pt]star.east)--node[above,scale=0.6]{$\scriptstyle{p_1^*}$}([yshift=2pt]ver13.west);
                 \draw[->] (ver13)--node[right]{$\scriptstyle{\lambda''}$}(X);
                 \draw[->] ([yshift=-2pt]Y.east)--([yshift=-2pt]X.west);
                 \draw[->] ([yshift=2pt]Y.east)--([yshift=2pt]X.west);

                 \draw[->] (star)--node[right]{$\scriptstyle{\lambda'}$}(Y);
                 \draw[->] (Z)--(Y);
                 \draw[->] (F)--node[right]{$\scriptstyle{\lambda}$}(Z);
                 \draw[->]
                 ([yshift=3.5pt]X.east)--([yshift=3.5pt]new.west);
                 \draw[->] ([yshift=-3.5pt]X.east)--([yshift=-3.5pt]new.west);
                 \draw[->] (X)--(new);
                 \draw[->]
                 ([yshift=3.5pt]ver13.east)--node[above,scale=0.3]{$\scriptstyle{p_{12}^*}$}([yshift=3.5pt]ver14.west);
                 \draw[->]
                 ([yshift=-3.5pt]ver13.east)--node[below,scale=0.3]{$\scriptstyle{p_{13}^*}$}([yshift=-3.5pt]ver14.west);
                 \draw[->]
                 (ver13)--node[above,scale=0.3]{$\scriptstyle{p_{23}^*}$}(ver14);
                 \draw[->]
                 (ver14)--node[right]{$\scriptstyle{\lambda'''}$}(new);
            \end{tikzpicture} 
        \end{equation}
        with exact rows — meaning
        that, for instance, the first row is exact iff $u^*$ induces
        an equivalence between $A(\sT)$ and the category of descent data in
        $A(\sT')$. Such descent data consists of pairs 
        $(T,\xi)$,
        where $T\in A(\sT')$ and $\xi\colon p_1^*(T)\xrightarrow{\cong}
        p_2^*(T)$ satisfying the cocycle condition:
        $(p_{23}^*\xi)\circ(p_{12}^*\xi)=p_{13}^*\xi$. Then
        $\lambda'$ fully faithful (resp. an equivalence) and $\lambda''$ faithful 
        (resp. fully faithful and $\lambda'''$ faithful) imply that $\lambda$
        is fully faithful (resp. an equivalence).

\textbf{Step 1, putting the colimit inside.}
 Let $Y^{\sh}_{\bar{x}}$ denote the strict henselization of $Y$ at
 $\bar{x}$. We first show that the natural functor
    \begin{equation}\label{étale nbhd}
        \varinjlim_{(U,\bar{u})}\FEt(\sX\times_{\F_q}U/M)\xrightarrow{\quad
        \simeq\quad}\FEt(\sX\times_{\F_q}Y^{\sh}_{\bar{x}}/M)
    \end{equation} 
    is an equivalence.

    For each $i$, choose a smooth surjective equivariant atlas
    $X_i\to\sX_i$ with $X_i$ a qcqs scheme; this is possible because
    $\sX_i$ is quasi-compact, and the Frobenius preserves every open subset
    of $X_i$. Put $X=X_1\times\cdots\times X_n$. Since $\sX$ is
    quasi-separated, the terms $X$, $X\times_{\sX}X$, and
    $X\times_{\sX}X\times_{\sX}X$ of the resulting \v{C}ech presentation
    are qcqs algebraic spaces.

    Let $T$ be any of these three terms. The limit results
    \cite[\href{https://stacks.math.columbia.edu/tag/07SK}{07SK},
    \href{https://stacks.math.columbia.edu/tag/084Z}{084Z},
    \href{https://stacks.math.columbia.edu/tag/07SL}{07SL}]{stacks-project}
    apply to the qcqs inverse system $T\times_{\F_q}U$ and give
    \[
       \varinjlim_{(U,\bar u)}\FEt(T\times_{\F_q}U)
       \simeq \FEt(T\times_{\F_q}Y^{\sh}_{\bar x}).
    \]
    The same assertion holds with the $M$-actions: the monoid
    $M\cong\N^n$ has finitely many generators and finitely many commutation
    relations, so the linearizations, their inverses, and all relations
    descend after increasing the index once.

    Finite \'{e}tale morphisms, with or without the $M$-linearization,
    satisfy smooth descent. The corresponding descent diagrams use only the
    three \v{C}ech terms displayed above. Moreover, a descent datum in the
    filtered colimit is represented, together with its finitely many cocycle
    identities, at one common stage. Applying the principle of Step~0 to this
    finite descent diagram proves \eqref{étale nbhd}.

    \textbf{Step 2, reduction to $\sX_i=X_i$ a connected affine $\F_q$-scheme of
    finite type.} To prove that \eqref{existence theorem equation} is an
    equivalence, set $A(-)\coloneqq
    \FEt(-\times_{\F_q}Y^{\sh}_{\bar x}/M)$ and
        $B(-)\coloneqq\FEt(-\times_{\F_q}\bar{x}/M)$. Choose atlases $u_i\colon
        X_i\to\sX_i$ and set $u\coloneqq u_1\times\dots\times u_n\colon
        (X,\rho_X)\to (\sX,\rho_{\sX})$, where $X\coloneqq X_1\times\cdots\times X_n$
        and $\rho_X$ is the partial Frobenius action. Each relevant \v{C}ech
        term is a finite product of qcqs algebraic spaces. After choosing
        finite affine étale covers of these factors and applying Step~0 to
        the resulting finite descent diagrams, we are reduced to the case
        where each $\sX_i=X_i$ is an affine scheme.
    Next, write each affine scheme $X_i$ as a cofiltered limit of affine schemes
    of finite type over $\F_q$, and reduce to the case where $X_i$ is of finite type over
    $\F_q$ by the limit results cited in Step~1. A finite-type affine scheme
    has only finitely many connected components, each preserved by Frobenius;
    applying Step~0 once more reduces us to a connected component.

    \textbf{Step 3, reduction to $Y$ excellent, reduced, and strictly
    henselian.} By Corollary \ref{Permanence v2}, the target category is
    unchanged if the field of the geometric point is replaced by an
    algebraic closure of the residue field at its image. We make this
    replacement at the outset.

    Since the problem is Zariski-local around $\bar{x}$, we may
    assume that $Y$ is affine. Writing $Y$ as a cofiltered limit of affine finite-type 
    $\F_q$-schemes
    and applying the limit results cited in Step~1,
    we may assume that $Y$ is of finite type over $\F_q$. Replacing $Y$ with
    its reduced closed subscheme, we may assume that $Y$ is reduced. Now
    replacing $Y$ by its strict henselization at $\bar{x}$ using \eqref{étale
    nbhd}, we may assume that $Y=\Spec(A)$, where $A$ is strictly henselian
    and reduced with separably closed residue field $\kappa$. The extension
    $k/\kappa$ is purely inseparable, hence
    $\Spec(k)\to\Spec(\kappa)$ is a universal homeomorphism. Invariance of
    the \'{e}tale site under universal homeomorphisms, together with Lemma
    \ref{universal homeomorphism vs cartesian}, allows us to replace
    $\bar{x}$ by $\Spec(\kappa)$. By
    \cite[\href{https://stacks.math.columbia.edu/tag/06LJ}{06LJ}]{stacks-project},
    $A$ is Noetherian; by
    \cite[\href{https://stacks.math.columbia.edu/tag/07QR}{07QR}]{stacks-project},
    $A$ is a G-ring; by \cite[\nopp 18.8.17, 18.7.5.1]{EGAIV4}, $A$ is
    universally catenary; and by \cite[\nopp 5.3]{Gre76} it is J-2. In
    conclusion, $A$ is excellent.

    \textbf{Step 4, the functor $\lambda$ is faithful and essentially surjective.}
Consider the following $2$-commutative diagram, where $\sX=X$ is affine and of finite type.
    \[
       \begin{tikzpicture}[xscale=6.9,yscale=-2.2]
            \node (A0_0) at (0, 0.5) {$\FEt(X/M_0)$};
            \node (A0_1) at (1, 0) {$\FEt(X\times_{\F_q}\Spec(A)/M)$};
            \node (A1_1) at (1, 1) {$\FEt(X\times_{\F_q}\Spec(\kappa)/M)$};
        
            \draw[->] (A0_0) -- node[anchor=south]{$\scriptstyle{w}$} (A0_1);
            
            \path (A0_1) edge[->] node[right] {$\scriptstyle{\lambda}$} (A1_1);
            \path (A0_0) edge[->] node[anchor=north west]{$\scriptstyle{\simeq}$} node[anchor=south west]
            {$\scriptstyle{v}$} (A1_1);
        \end{tikzpicture}
    \]
    We must show that $\lambda$ is an equivalence. By Remark
    \ref{dropping out any i}, the category in the lower right may be described
    using the action of $M_0\times\phi_\kappa^\N$. Choose an algebraic closure
    $\bar\kappa$ of $\kappa$. Since $\kappa$ is separably closed,
    $\bar\kappa/\kappa$ is purely inseparable. Invariance under universal
    homeomorphisms and Corollary \ref{Permanence v2} give
    \[
       \FEt(X\times_{\F_q}\Spec(\kappa)/M)
       \simeq \FEt(X\times_{\F_q}\Spec(\bar\kappa)/M)
       \simeq \FEt(X/M_0).
    \]
    Therefore $v$ is an equivalence. By Lemma \ref{connectedness} and Lemma
    \ref{Galois category}, all the categories in the above diagram are Galois
    categories; since all the functors are maps of Galois categories
    (i.e., exact), they are faithful.
    It follows that $w$ is fully faithful and $\lambda$ is essentially
    surjective. Thus, to conclude, it suffices to show that $\lambda$ is fully
    faithful, or equivalently that $w$ is essentially surjective.

    \textbf{Step 5, reduction to $X_i$ and $Y$ normal, with $Y$ having an algebraically
    closed field of fractions.} The normalization morphisms used below are
    finite because all rings at this stage are excellent. They are therefore
    proper, surjective, and of finite presentation. Such morphisms are of
    effective descent for finite \'{e}tale covers by
    \cite[Expos\'{e} IX, Th\'{e}or\`{e}me 4.12]{SGA1}; the same is true for
    fpqc morphisms. We shall use this descent without further comment.
    
    As each $X_i$ is excellent, its
    normalization map is finite. Apply diagram \eqref{descent on X1} to
    $A(-),B(-)$ as in Step 2, with $u$ the
    product of normalization maps. After decomposing the normalizations and
    their \v{C}ech products into their finitely many $G_0$-components, note that
    $(X'',\rho_{X''})=\bigsqcup_{i=1}^r(Z_i,\rho_{Z_i})$ is a finite disjoint union of
    $G_0$-components, so $\lambda''$ is a product of functors
    \(
    \FEt(Z_i\times_{\F_q}Y/M)\longrightarrow\FEt(Z_i\times_{\F_q}\bar{x}/M),
    \)
    each of which is faithful (as the source and target are both Galois categories, cf.~ Lemma \ref{connectedness} and Theorem
    \ref{Galois category}). Therefore,
    $\lambda''$ is faithful, and we may assume that each $X_j$ is normal. 

    For $Y$, we first note that there are finitely many irreducible components,
    each of which is the spectrum of a strictly henselian excellent domain.
    Indeed, if $Y_i=\Spec(A/\mathfrak p_i)$, then $A/\mathfrak p_i$ is a local
    quotient of the henselian ring $A$, has the same separably closed residue
    field, and is excellent.
    If $Y$ has $m$ irreducible components $Y_1,\ldots,Y_m$, set
    $Y'\coloneqq Y_1\sqcup (Y_2\cup\cdots\cup Y_m)$, which is an h-covering of
    $Y$. Consider the following diagram (diagram \eqref{descent
    on X1} ``with the other factor'')
\begin{equation}\label{descent on X2}
                    \begin{tikzpicture}[xscale=3.5,yscale=1.2,baseline=(current  bounding  box.center)]
                 \path
                 (0,0)     node[scale=1] (F) {$\FEt(X\times_{\F_q}Y/M)$}
                     ++(0:1.5)  node[scale=1] (star)
                     {$\FEt(X\times_{\F_q}Y'/M)$}
                     ++(0:1.5) node[scale=1] (ver13)
                     {$\FEt(X\times_{\F_q}Y''/M)$}
                     ++(-90:1.5) node[scale=1] (X)
                     {$\displaystyle\FEt(X\times_{\F_q}\bar{x}''/M)$}
                     ++(-180:1.5)  node[scale=1] (Y)
                     {$\displaystyle\FEt(X\times_{\F_q}\bar{x}'/M)$}
                     ++(-180:1.5) node[scale=1] (Z)
                     {$\displaystyle\FEt(X\times_{\F_q}\bar{x}/M)$};
                 \draw[->] (F)--(star);
                 \draw[->] ([yshift=-3pt]star.east)--([yshift=-3pt]ver13.west);
\draw[->] ([yshift=3pt]star.east)--([yshift=3pt]ver13.west);
                 \draw[->] (ver13)--node[right]{$\scriptstyle{\lambda''}$}(X);
                 \draw[->] ([yshift=-3pt]Y.east)--([yshift=-3pt]X.west);
                 \draw[->] ([yshift=3pt]Y.east)--([yshift=3pt]X.west);

                 \draw[->] (star)--node[right]{$\scriptstyle{\lambda'}$}(Y);
                 \draw[->] (Z)--(Y);
                 \draw[->] (F)--node[right]{$\scriptstyle{\lambda}$}(Z);         
            \end{tikzpicture} 
        \end{equation}
        where $\bar{x}'\coloneqq \bar{x}\sqcup\bar{x}$ and $\bar{x}''\coloneqq
        \bar{x}'\times_{\bar{x}}\bar{x}'$. By the effective descent theorem
        cited at the beginning of this step (and its base change to $X$), the
        rows in \eqref{descent on X2} are exact. The finitely many
        $M$-linearizations and their compatibility relations descend as well.
        Note that $\lambda''$ is faithful, as
        the residue map $\bar{x}''\to Y''$ induces a bijection on connected
        components ($Y''$ is a finite disjoint union of spectra of strictly henselian local
        rings by \cite[\href{https://stacks.math.columbia.edu/tag/04gh}{04GH}]{stacks-project}).
        Using induction on $m$, we are reduced to the case
    where $A$ is a domain.

    Now let $K$ be the field of fractions of $A$, and let $B$ be the
    normalization of $A$ in some finite extension $L/K$. Since $A$ is
    excellent, $B$ is finite over $A$. Moreover, $B$ and
    $C\coloneqq B\otimes_A B$ are local and strictly henselian. Indeed,
    \cite[\href{https://stacks.math.columbia.edu/tag/04GH}{04GH}]{stacks-project}
    writes the finite $A$-algebra $B$ as a product of henselian local rings;
    since $B$ is a domain, there is only one factor. Put
    $D=B/\mathfrak m_A B$, let $\mathfrak n$ be its maximal ideal, and write
    $\kappa_B$ for its residue field. Then $\kappa_B/\kappa$, where
    $\kappa=A/\mathfrak m_A$, is finite purely inseparable. The ideal
    $\mathfrak n\otimes_\kappa D+D\otimes_\kappa\mathfrak n$ is nilpotent in
    $D\otimes_\kappa D$, and the quotient is
    $\kappa_B\otimes_\kappa\kappa_B$, which has a unique prime ideal.
    Therefore $D\otimes_\kappa D$ is Artinian local, and $C$ has exactly one
    prime ideal lying over $\mathfrak{m}_A$. Every maximal ideal of the finite
    $A$-algebra $C$ lies over $\mathfrak m_A$, so $C$ is local. Since $C$ is finite over the
    henselian local ring $B$, it is henselian by
    \cite[\href{https://stacks.math.columbia.edu/tag/04GH}{04GH}]{stacks-project}.
    The residue field $\kappa_C$ of $C$, like $\kappa_B$, is a finite purely
    inseparable extension of the separably closed field $\kappa$, and hence
    is separably closed.
    Thus both $B$ and $C$ are strictly henselian.

    Fix an algebraic closure $\Omega$ of $\kappa$. The $\kappa$-embeddings
    $\kappa_B\to\Omega$ and $\kappa_C\to\Omega$ are unique; in particular,
    they are compatible with both maps $B\rightrightarrows C$. Hence
    $\Spec(\Omega)$ gives compatible geometric closed points of
    $\Spec(A)$, $\Spec(B)$, and $\Spec(C)$. Since $\Omega/\kappa$ is purely
    inseparable, invariance under universal homeomorphisms shows that replacing
    the closed fiber over $\kappa$ by the one over $\Omega$ does not change
    the relevant finite étale $M$-categories. Moreover, the restriction
    functor
    \[
       \FEt(X\times_{\F_q}\Spec(C)/M)\longrightarrow
       \FEt(X\times_{\F_q}\Spec(\Omega)/M)
    \]
    is faithful: componentwise, its source and target are Galois categories
    by Corollary \ref{connectedness} and Theorem \ref{Galois category}, and
    the restriction functor is exact. Applying Step~0 to the normalization
    map $\Spec(B)\to\Spec(A)$, with these compatible geometric points, reduces
    us to the case where $Y$ is normal.

    Let $\bar{Y}=\Spec(\bar{A})$ be the normalization of $Y$ in an algebraic
    closure $\bar{K}$ of $K$. The ring $\bar A$ is the filtered colimit of the
    finite normalizations $B$ inside the finite subextensions of $\bar K/K$.
    These are strictly henselian by the preceding argument, and the transition
    maps are local. The ring $\bar A$ is normal by construction and strictly
    henselian by
    \cite[\href{https://stacks.math.columbia.edu/tag/04GI}{04GI}]{stacks-project}.
    The unique compatible embeddings of their residue fields into $\Omega$
    also define a geometric closed point $\Spec(\Omega)\to\Spec(\bar A)$.
    Suppose $T_1,T_2\in\FEt(X\times_{\F_q}Y/M)$ and $t\colon
    T_1\otimes_A\kappa\to T_2\otimes_A\kappa$ is a map in
    $\FEt(X\times_{\F_q}\bar{x}/M)$. Suppose that $t$ lifts to a morphism
    $\bar a$ over $\bar A$. The morphism and its finitely many
    generator-equivariance identities are defined over some finite
    normalization $A\subseteq B\subseteq\bar A$. Its two pullbacks to
    $C=B\otimes_A B$ have the same restriction to $\Spec(\Omega)$, since both
    restrict to $t_\Omega$. By the faithfulness just proved for $C$, they are
    equal. Effective descent along $\Spec(B)\to\Spec(A)$ therefore yields a
    unique morphism $a$ over $A$ lifting $t$.
    Thus we are reduced to the case where the fraction field $K$ of $A$ is
    algebraically closed and $A$ is normal (but not necessarily Noetherian).

    \textbf{Step 6, final conclusion.} Recall that for any integral $\F_q$-algebra
    $R$, the integral closure $R_{\et}$ of $\F_q$ inside $R$ is a possibly
    infinite algebraic extension of $\F_q$, hence a filtered colimit of
    finite \'{e}tale $\F_q$-algebras. The scheme $\Spec(R)$ is geometrically connected over $R_{\et}$
    (cf.~\cite[Lem. 2.7]{TZ1}). If $R$ is normal with field of fractions $L$, then
    $R_{\et}=L_{\et}$. Apply this with $R=A$ and $L=K$, and put
    $k_0=A_{\et}=K_{\et}$. Since $K$ is algebraically closed, so is $k_0$:
    a root in $K$ of a polynomial over $k_0$ is algebraic over $\F_q$ and
    hence belongs to $K_{\et}$. Both $\Spec(A)$ and $\Spec(K)$ are geometrically
    connected over $k_0$. The finitely many connected components of
    $X\times_{\F_q}\Spec(k_0)$ are geometrically connected, because $k_0$ is
    algebraically closed. Their products with $\Spec(A)$ and $\Spec(K)$ remain
    connected. Thus base change gives bijections on connected components, and
    the schemes $X\times_{\F_q}\Spec(A)$ and
    $X\times_{\F_q}\Spec(K)$ have the same number of connected components.
    After the preceding normalization, every $X_i$ is a normal finite-type
    scheme over the perfect field $\F_q$, hence geometrically normal by
    \cite[\href{https://stacks.math.columbia.edu/tag/038O}{038O}]{stacks-project}.
    Iterating
    \cite[\href{https://stacks.math.columbia.edu/tag/06DG}{06DG}]{stacks-project}
    after arbitrary field extension shows that $X$ is geometrically normal;
    another application shows that $X\times_{\F_q}\Spec(A)$ is normal. For each
    connected component $Z$ of this scheme, the preceding component description
    and normality show that $Z$ is normal and integral. Flatness and going down
    show that $Z$ meets the generic fiber, and equality of the component counts
    shows that $Z_K$ is connected; being normal, it is integral. A morphism
    between finite \'{e}tale covers over $Z_K$ restricts to one over the common
    function field $k(Z)$. By
    \cite[\href{https://stacks.math.columbia.edu/tag/0BQM}{0BQM}]{stacks-project}
    it extends uniquely over $Z$, and its restriction back to $Z_K$ is the
    original morphism by faithfulness of restriction from the normal integral
    scheme $Z_K$ to $k(Z)$. Componentwise, this shows that
    the functor
    \(
    \FEt(X\times_{\F_q}\Spec(A))\to\FEt(X\times_{\F_q}\Spec(K))
    \)
    is fully faithful. It follows that
    \(
    \FEt(X\times_{\F_q}\Spec(A)/M)\xrightarrow{\ \theta\ }\FEt(X\times_{\F_q}\Spec(K)/M)
    \)
    is also fully faithful by the rigidification observation following Lemma
    \ref{universal homeomorphism vs cartesian}. Since
    \(
    \FEt(X/M_0)\xrightarrow{\ \delta\ }\FEt(X\times_{\F_q}\Spec(K)/M)
    \)
    is an equivalence by Remark \ref{dropping out any i} and Corollary
    \ref{Permanence v2},
    the equality $\theta\circ w=\delta$ implies that $w$ — which is
    already fully faithful by Step 4 — is an equivalence, as desired.
\end{proof}

\begin{rmk}\label{non-qc counterexample}
The quasi-compactness hypothesis in Lemma \ref{pi_1 existence theorem} cannot
be omitted, even for $n=1$ and schemes. Take $M=\phi_X^\N$, and let $S$ be
the spectrum of the ring of eventually constant
$\F_q$-valued sequences. Its points are the isolated points
$1,2,\ldots$ and a point $\infty$, every neighborhood of which contains all
but finitely many isolated points. Let $C$ be a split nodal rational curve,
$V\subset C$ the complement of its node, and let $D\to C$ be the geometrically connected
double cover obtained by cross-gluing two copies of the normalization; its
restriction to $V$ is the trivial double cover. Glue copies $C_i$ along their
common open $V$ to obtain a connected non-quasi-compact scheme $X$.

On $C_i\times S$, use $D$ over the clopen point $i$ and the trivial double
cover over its complement. The resulting covers glue over $V\times S$ to a
finite \'{e}tale double cover $E\to X\times S$. Put
$W\coloneqq(X\sqcup X)\times_{\F_q}S$. The $q$-Frobenius of $S$ is the
identity, so the absolute Frobenius of $X\times S$ is
$\phi_X\times\id_S$. Thus $E$ and $W$ carry their canonical
$M$-linearizations. The fiber of $E$ over a geometric point above $\infty$
is the trivial cover $W_\infty$, whereas the fiber $E_i$ over every isolated
geometric point $i$ is connected. The
sheetwise $M$-equivariant isomorphism $E_\infty\simeq W_\infty$ cannot spread
to an \'{e}tale neighborhood of $\infty$. Thus the functor in
\eqref{existence theorem equation} is not full. This is precisely the failure
excluded by the finite qcqs presentation in Step~1.
\end{rmk}

\begin{lem}[Stein factorization] Suppose that $\sX_1,\ldots,\sX_n$ are qcqs
    and that $\sX$ is $G_0$-connected, and let $\sX_{n+1}$ be an arbitrary
    algebraic stack. For any
    $(\sY,\rho_{\sY})\in\FEt(\sX\times_{\F_q}\sX_{n+1}/M)$, there
    exists $\sT\in\FEt(\sX_{n+1})$ and a 2-commutative diagram 
    \begin{equation}\label{stein factorization}
        \begin{tikzpicture}[xscale=2.0,yscale=1.2,bmr/.pic={\draw (0,0)--++(-90:2mm)--++(180:2mm);},baseline=(current  bounding  box.center)]
             \path
                    (0,0)     node (F) {$\sY$}
                 %+(-43:.4) pic[scale=1,red]{bmr}
                 +(0:1.5)  node (star) {$\sT$}
                 ++(-90:1.5) node (X) {$\sX\times_{\F_q}\sX_{n+1}$}
                 +(0:1.5)  node (Y) {$\sX_{n+1}$};
             \draw[->] (F)--(star);
             \draw[->] (F)--(X);
             \draw[->] (X)--(Y) node[midway,above,scale=.6]{};
             \draw[->] (star)--(Y);
        \end{tikzpicture} 
    \end{equation}
     where $\sY\to \sT$ is $M$-equivariant
     (with the trivial $M$-action on $\sT$), and each geometric
     fiber $\sY_{\bar{t}}$ of $\sY\to \sT$ is $G$-connected. If
     diagram \eqref{stein factorization} exists, then it satisfies the following 
     universal property: For any $\sT'\in\FEt(\sX_{n+1})$ and any $M$-equivariant
     map $\sY\to \sT'$ forming a diagram analogous to \eqref{stein factorization}, 
     there is a unique arrow $\lambda\colon\sT\to \sT'$ making the two diagrams 
     compatible.
\end{lem}

\newcounter{Stein}
\setcounter{Stein}{1}

\begin{proof} 
    The general strategy is as follows. We first show that if $\sT$ exists, it must satisfy 
    the universal property and is therefore unique up to a unique isomorphism. 
    By fppf descent for finite étale covers, it suffices to prove existence fppf locally, 
    which allows us to make various reductions.

    \textbf{Step \theStein\stepcounter{Stein}: Universal property.} 
    Assume that $\sT$ exists. We prove its universal property: For any $M$-equivariant 
    map $\sY\to \sT'$ over $\sX_{n+1}$, there exists a unique map $\lambda\colon\sT\to \sT'$
    making the diagram \eqref{stein factorization} for $\sT$ compatible with that for $\sT'$.
    
    Since $\sY\to \sT$ is faithfully flat, such a map $\lambda$ is unique if it exists 
    (note that $\FEt(\sX_{n+1})$ is a $1$-category, so this makes sense). 
    By fppf descent, we may construct $\lambda$ locally on $\sX_{n+1}$. 
    Because every finite étale cover splits étale locally, we may assume that 
    $\sT = \bigsqcup_{i=1}^r \sX_{n+1}$ is a finite disjoint union. 
    Working with each component separately, we may further assume $\sT = \sX_{n+1}$; 
    consequently, each geometric fiber of $\sY \to \sX_{n+1}$ is $G$-connected.
    
    The map $f\colon \sY \to \sT'$ is universally open. Indeed, the morphism
    $(p_{\sX},f)\colon\sY\to\sX\times_{\F_q}\sT'$ is a morphism between
    finite étale covers of $\sX\times_{\F_q}\sX_{n+1}$, hence is finite
    étale. The projection $\sX\times_{\F_q}\sT'\to\sT'$ is universally open:
    this may be checked after pulling back to a smooth atlas of $\sX$, where it
    follows from
    \cite[\href{https://stacks.math.columbia.edu/tag/0383}{0383}]{stacks-project}.

    Let $\sU\subseteq\sT'$ be the open image of $f$. For every geometric point
    $\bar s$ of $\sX_{n+1}$, the finite discrete set $\sT'_{\bar s}$ has
    trivial $G$-action. Since $\sY_{\bar s}$ is nonempty and $G$-connected,
    its equivariant image in $\sT'_{\bar s}$ is a singleton. Consequently,
    the étale morphism $u\colon\sU\to\sX_{n+1}$ is surjective and radicial,
    hence an isomorphism. The required map $\sT=\sX_{n+1}\to\sT'$ is the
    composite of $u^{-1}$ with the open immersion $\sU\hookrightarrow\sT'$.

    \textbf{Step \theStein\addtocounter{Stein}{-1}: Reduction to $\sX_{n+1}$ affine.} 
    By Step \theStein\addtocounter{Stein}{2}, it suffices to prove existence fppf locally. 
    Thus we may assume that $\sX_{n+1} = X_{n+1} = \Spec(A)$ is affine.

    \textbf{Step \theStein\addtocounter{Stein}{1}: Reduction to $\sY$ being a pullback from $\FEt(\sX/M_0)$.} 
    By Lemma \ref{pi_1 existence theorem}, Remark \ref{dropping out any i},
    and Corollary \ref{Permanence v2}, for any geometric point
    $\bar{x}$ of $X_{n+1}$, we can find an affine étale neighborhood $U \to X_{n+1}$ such that 
    the pullback $(\sY_U,\rho_{\sY_U}) \in \FEt(\sX\times_{\F_q}U/M)$ of $(\sY,\rho_{\sY})$ 
    comes from an object $\sZ \in \FEt(\sX/M_0)$. 
    These neighborhoods form an étale covering. By the universal property in
    Step~1, the resulting local factorizations have unique identifications on
    overlaps, and these identifications satisfy the cocycle condition. Hence
    fppf descent glues them, so it is enough to construct $\sT$ over $U$.

    \textbf{Step \theStein\stepcounter{Stein}: Final conclusion.} 
    Let $\sZ_1,\sZ_2,\dots,\sZ_m$ be the $G_0$-components of $\sZ$ 
    (recall that $\FEt(\sX/M_0)$ is Galois). 
    Set $\sT \coloneqq \bigsqcup_{i=1}^m U$, and define the map $\sY \to \sT$ as the projection 
    $\delta\colon \bigsqcup_{i=1}^m (\sZ_i \times_{\F_q} U) \to \bigsqcup_{i=1}^m U$, 
    where we use the identification 
    \[
    \FEt(\sX\times_{\F_q}U/M) \simeq \FEt(\sX\times_{\F_q}U/(M_0 \times \phi_U^{\N}))
    \]
    to ensure that the projection is equivariant with the trivial action on $U$.
    
    For a geometric point $\bar t=(\bar x,i)$ of the $i$-th copy of $U$ in
    $\sT=\bigsqcup_iU$, the pullback of
    $(\sZ_i,\rho_{\sZ_i})$ to $\FEt(\sX\times_{\F_q}k/M)$ is $G$-connected, as by  Corollary \ref{Permanence v2}
    \[
    \FEt(\sX\times_{\F_q}k/M) \simeq \FEt(\sX\times_{\F_q}k/(M_0 \times \phi_k^{\N}))
    \simeq \FEt(\sX/M_0)
    \]
    and $(\sZ_i,\rho_{\sZ_i})$ is $G_0$-connected. Hence $\delta^{-1}(\bar t)$ is $G$-connected,
    as required.
\end{proof}

\begin{thm}\label{exact sequence thm} Suppose that $\sX_i$ $(1\leq i\leq n)$
    are connected qcqs algebraic stacks and that $\sX_{n+1}$ is an arbitrary
    connected algebraic stack over $\F_q$. Then for any geometric point
    $\bar{x}$ of $\sX\times_{\F_q}\sX_{n+1}$, the top sequence in
      \begin{equation}\label{exact sequence}
             \begin{tikzpicture}[xscale=4.8,yscale=-2.2,baseline=(current  bounding  box.center)]
                  \node (A0_-1) at (-.6, 0) {$1$};
                  \node (A0_0) at (0, 0) {$\pi_1^\et(\sX\times_{\F_q}k/M,\bar{x})$};
                  \node (A0_1) at (1, 0) {$\pi_1^\et(\sX\times_{\F_q}\sX_{n+1}/M,\bar{x})$};
                  \node (A0_2) at (1.9, 0) 
                  {$\pi_1^\et(\sX_{n+1},\bar{x})$};
                  \node (A1_0) at (0, 1)
                  {$\pi_1^\et(\sX/M_0,\bar{x})$};
                  \node (A0_3) at (2.4, 0) {$1$};
                  %\node (A1_1) at (1, 1) {$right-lower$};
              \draw[->] (A0_0) -- node[anchor=south]{$\scriptstyle{\gamma}$} (A0_1);
                  \draw[->] (A0_-1) --
                  node{$\scriptstyle{}$} (A0_0);
                  \draw[->] (A0_2) --
                  node[above]{$\scriptstyle{}$} (A0_3);
                  \draw[->] (A0_1) -- node[below]{$\scriptstyle{\alpha}$} (A1_0);
          
                  \path (A0_1) edge[->] node[above]
                  {$\scriptstyle{\beta}$} (A0_2);
                  \path (A0_0) edge[->] node[left]
                  {$\scriptstyle{\cong}$} (A1_0);
              \end{tikzpicture}
                \end{equation}
      is exact, where $k$ is the field of definition of $\bar{x}$.
\end{thm}
\begin{proof} Thanks to Corollary \ref{Permanence v2}, the natural map
    $\alpha$ is a retraction of $\gamma$, so that $\gamma$ is injective. 
    Since Corollary
    \ref{connectedness} assures that the pullback of a connected cover in $\FEt(\sX_{n+1})$ to
    $\FEt(\sX\times_{\F_q}\sX_{n+1}/M,\bar{x})$ is $G$-connected, by
    \cite[\href{https://stacks.math.columbia.edu/tag/0BN6}{0BN6}]{stacks-project}
    $\beta$ is surjective.
    We are left to show the middle
    exactness.
    By
    \cite[\href{https://stacks.math.columbia.edu/tag/0BTQ}{Lemmas 58.4.3 and 58.4.5}]{stacks-project},
    it is enough to verify the following connected-object condition. Lemma
    58.4.5 then gives normality of $\operatorname{im}(\gamma)$, while Lemma
    58.4.3 identifies the objects inflated from the base and yields
    $\operatorname{im}(\gamma)=\ker(\beta)$.

    \texttt{
    Let $(\sY,\rho_{\sY})\in\FEt(\sX\times_{\F_q}\sX_{n+1}/M)$ be a
    connected object. Suppose the pullback
    $(\sY_{\bar{x}},\rho_{\sY_{\bar{x}}})$ of
    $(\sY,\rho_{\sY})$ to $\FEt(\sX\times_{\F_q}k/M)$ has a
    $G$-component which maps
    isomorphically to $\sX\times_{\F_q}k$, then the Stein
    factorization of $(\sY,\rho_{\sY})$ induces a Cartesian diagram. 
    \[
        \begin{tikzpicture}[xscale=2.0,yscale=1.2,bmr/.pic={\draw (0,0)--++(-90:2mm)--++(180:2mm);}]
             \path
                    (0,0)     node (F) {$\sY$}
                 +(-43:.4) pic[scale=1,red]{bmr}
                 +(0:1.5)  node (star) {$\sT$}
                 ++(-90:1.5) node (X) {$\sX\times_{\F_q}\sX_{n+1}$}
                 +(0:1.5)  node (Y) {$\sX_{n+1}$};
             \draw[->] (F)--(star);
             \draw[->] (F)--(X);
             \draw[->] (X)--(Y) node[midway,above,scale=.6]{};
             \draw[->] (star)--(Y);
        \end{tikzpicture}
    \]
}
    One has to show that $\lambda\colon \sY\to
    \sX\times_{\F_q}\sT$ is an isomorphism. For each geometric point $\bar{s}$ of $\sT$ the
    pullback $\sX\times_{\F_q}\kappa(\bar{s})$ is $G$-connected by
    Corollary \ref{connectedness}. It follows
     that
 $\lambda_{\bar{s}}\colon \sY_{\bar{s}}\to \sX\times_{\F_q}\kappa(\bar{s})$ as well as
     $\lambda$ are surjective ($\sY_{\bar{s}}$ is a $G$-component, hence is not empty).
    The morphism $\lambda$ is a morphism between finite étale covers of
    $\sX\times_{\F_q}\sX_{n+1}$, hence is finite étale. Since it is
    surjective and equivariant, the $G$-connectedness of $\sY$ implies that
    $\sX\times_{\F_q}\sT$ is $G$-connected. The rank loci of $\lambda$ are
    $G$-stable and clopen, so $\lambda$ has constant degree.
    Therefore to show that $\lambda$ is an isomorphism, it is enough to show that
    the degree of $\lambda$ equals $1$.

    The $G$-component of $(\sY_{\bar{x}},\rho_{\sY_{\bar{x}}})$ which goes
    isomorphically to
    $\sX\times_{\F_q}k$ is mapped to a geometric point
    $\bar{s}\in\sT(k)$ lying over $\bar{x}\in\sX_{n+1}(k)$. It is moreover
    the inverse image of $\bar{s}$ by the very construction of the Stein
    factorization.  Thus 
    $\lambda$ is an isomorphism at $\bar{s}$, and the degree of $\lambda$
    is 1, as desired.
    \end{proof}

    \begin{proof}[Proof of Theorem \ref{Drinfeld's Lemma}] (1) is a special
        case of Theorem \ref{Galois category}, so let's consider (2). The case
        $n=1$ follows from Lemma \ref{are mutually commute}, since then
        $M_0$ is the trivial monoid. For $n\geq2$, after
        permuting the factors, we may suppose that
        $\sX_1,\ldots,\sX_{n-1}$ are qcqs; the last factor is allowed to be
        non-qcqs.
        By Remark \ref{dropping out any i}, Theorem \ref{exact sequence thm},
        and diagram \eqref{exact sequence}
        we have
        \[\pi_1^\et(\sX/\phi_1^\N\times\cdots\times\phi_{n-1}^\N)\xrightarrow{\quad\cong\quad}\pi_1^\et(\sX_1\times\cdots\times\sX_{n-1}/\phi_1^\N\times\cdots\times\phi_{n-2}^\N)\times\pi_1^\et(\sX_n).\]
        Here $\alpha\circ\gamma=\id$ and
        $\ker(\beta)=\operatorname{im}(\gamma)$, so $(\alpha,\beta)$ is an
        isomorphism onto the displayed direct product. Induction on $n$ now
        gives the result.
    \end{proof}

%\sloppy
\printbibliography

\end{document}